%% file: MF-ArXiv-2024-01-03.tex
\documentclass[review,onefignum,onetabnum]{siamonline220329}

\input{MF_shared}

\ifpdf
\hypersetup{
  pdftitle={Interspecific density-dependent model},
  pdfauthor={MTAR, FEKIH}
}
\fi

\begin{document}

\maketitle
\begin{abstract}
This paper studies a two microbial species model in competition for a single resource in the chemostat including general interspecific density-dependent growth rates with distinct removal rates for each species.
We give the necessary and sufficient conditions of existence, uniqueness, and local stability of all steady states.
We show that a positive steady state, if it exists, then it is unique and unstable.
In this case, the system exhibits a bi-stability where the behavior of the process depends on the initial condition.
Our mathematical analysis proves that at most one species can survive which confirms the competitive exclusion principle.
We conclude that adding only interspecific competition in the classical chemostat model is not sufficient to show the coexistence of two species even considering mortality in the dynamics of two species.
Otherwise, we focus on the study, theoretically and numerically, of the operating diagram which depicts the existence and the stability of each steady state according to the two operating parameters of the process which are the dilution rate and the input concentration of the substrate.
Using our mathematical analysis, we construct analytically the operating diagram by plotting the curves that separate their various regions.
Our numerical method using MATCONT software validates these theoretical results but it reveals new bifurcations that occur by varying two parameters as Bogdanov-Takens and Zero-Hopf bifurcations.
The bifurcation analysis shows that all steady states can appear or disappear only through transcritical bifurcations.
\end{abstract}
\begin{keywords}
Bi-stability, Bifurcation, Chemostat, MATCONT, Operating Diagram.
\end{keywords}
\begin{MSCcodes}
  34A34, 34D20, 37N25, 92B05
\end{MSCcodes}
\section{Introduction}	
Invented by Monod \cite{Monod1950} and almost simultaneously by Novick and Szilard \cite{Novick1950} in 1950, the chemostat is a bioreactor fed continuously with a nutrient flow for the cultivation of cells or microorganisms. More precisely, the chemostat represents a device consisting of an enclosure, containing the reaction volume, connected by input to feed the system with resources, and an outlet through which microorganisms are removed. 
It is the source of several mathematical models describing various biological and ecological interactions between microbial species in the chemostat, for instance: mutualism \cite{VandermeerJTB1978,VetSR2020}, commensalism \cite{BenyahiaJPC2012,BernardBB2001,SbarciogBE2010}, syntrophy \cite{DaoudMMNP2018,ElhajjiMBE2010,FekihSIADS2021,SariMB2016,XuJTB2011}, predator-prey relationship \cite{BoerMB1998,LiJMAA2000}, etc.
Some of these models have been studied to understand and exploit the anaerobic digestion process used in wastewater and waste treatment. It includes a large number of species that coexist with a very complex relationship (see \cite{FekihSIADS2021,MeadowsSIAP2019}, and the references therein).

The classic chemostat model describes the competition of two or more microbial populations on a single non-reproducing nutrient \cite{SmithBook1995}. 
It is well known that the \textit{Competitive Exclusion Principle} (CEP) predicts, in a generic case, that only one of them can survive. That is, only the species with the lowest "break-even" concentration survives while all other species will be excluded \cite{SmithBook1995}.
Although this theoretical prediction has been corroborated by the experiences of Hansen and Hubell \cite{HansenSc1980}, this principle (CEP) is not compatible with the great biodiversity found in nature.
Thus, to reconcile the theory and the experimental results, several suggestions and new hypotheses were considered in the classic chemostat model.
More precisely, to promote the coexistence of species in microbial ecosystems, many recent researches take into account various mechanisms, such as the inter and intraspecific competition \cite{AbdellatifMBE2016,DeLeenheerJMAA2006}, the flocculation of the species or wall attachment \cite{FekihJMAA2013,FekihAMM2016,FekihSIADS2019,HaegemanJBD2008,PilyuginSIAP1999}, the density-dependence of the growth functions \cite{BorsaliIJB2015,FekihMB2017,LobryCRBio2006EN,LobryEJDE2007,LobryCRMA2005,LobryCRBio2006,MtarIJB2021,MtarDcdsB2022}, etc.

Historically, an extension of the classical chemostat model was considered in a series of papers by Lobry et al. \cite{HarmandJPC2008,LobryCRBio2006EN,LobryEJDE2007,LobryCRMA2005,LobryCRBio2006}. 
This model is the following one
\begin{equation}                   \label{ModelLobry}
\left\{\begin{array}{lll}
\dot{S}   &=&  \ds D(S_{in}-S)- \sum_{i=1}^{n} \mu_i(S,X_1,\ldots,X_n)X_i,\\ 
\dot{X}_i &=&  (\mu_i(S,X_1,\ldots,X_n)-D_i)X_i, \quad i=1,\ldots,n
\end{array}\right.
\end{equation}
where $S$ and $X_i$ denote, respectively, the concentrations of the substrate and the concentration of species $i$ present in the reactor at time $t$; $S_{in}$ and $D$ denote, respectively, the concentration of substrate in the feed bottle and the dilution rate of the chemostat; $\mu_i(S,X_1,\ldots,X_n)$ represents the density-dependent growth rate of species $i$ which it depends not only on $S$ but also on the concentrations of species $i= 1,\ldots,n$.
The function $\mu_i$ is assumed to be increasing in the variable $S$ and decreasing in each variable $X_j$.
The removal rate $D_i$ can be interpreted as the sum of the dilution rate and the natural death rate of species $i$.
In \cite{LobryEJDE2007,LobryCRMA2005,LobryCRBio2006}, the authors considered the case where $\mu_i(S,X_i)$ depends only on $S$ and the concentration of species $i$ (that is, only intraspecific interference is present).
They proved that the coexistence of all species is possible.
A global study of model \cref{ModelLobry} with $D_i=D$ for all species was done in \cite{LobryCRMA2005}.
The case with distinct removal rates ($D_i\neq D$ for all species) was considered in \cite{LobryEJDE2007,LobryCRBio2006}.
For particular growth rates, general model \cref{ModelLobry} was studied only numerically in \cite{HarmandJPC2008,LobryCRBio2006EN}.
They demonstrated that the coexistence of species can occur for small enough interspecific interference but it is no longer the case when it is sufficiently large.
It was proved how the intraspecific interference is a mechanism of coexistence for competing species in the chemostat.

In this context, El Hajji et al. \cite{ElhajjiIJB2018,ElhajjiJBD2009} focused on the effect of interspecific interference on the coexistence of microbial species in the chemostat.
They considered a particular case of model \cref{ModelLobry} with $n=2$.
With the same dilution rate, the two-species model is written as follows:
\begin{equation}                              \label{ModelDEgaux}
\left\{\begin{array}{lll}
\dot{S}   &=&  D(S_{in}-S) - \mu_1(S,X_2)X_1-\mu_2(S,X_1)X_2, \\
\dot{X}_1 &=&  (\mu_1(S,X_2)-D)X_1,  \\ 
\dot{X}_2 &=&  (\mu_2(S,X_1)-D)X_2,
\end{array}\right.
\end{equation}
where the growth function $\mu_i(S,X_j)$ can depends now only on the substrate $S$ and the concentration of the species $X_j$, $i\neq j$.
In \cite{ElhajjiIJB2018}, three types of interspecific interactions between species are considered to explain the coexistence.
The first type is the predator-prey relationship which is characterized by the fact that one species (the prey) promotes the growth of the other species (the predator) which in turn inhibits the growth of the prey species.
In \cite{ElhajjiIJB2018}, it has shown that system \cref{ModelDEgaux} may exhibit the coexistence with a multiplicity of positive steady states.
In this case of interaction, our study in \cite{MtarDcdsB2022} has shown that adding the mortality of species in system \cref{ModelDEgaux} can destabilize the positive steady state through Hopf bifurcation where the coexistence can be around a stable limit cycle and not only around a positive steady state as the case without mortality.
The second type is the obligate mutualistic relationship where each species promotes the growth of the other species. That is, the function $\mu_i$ is increasing in $S$ and $X_j$.
Using general monotonic growth rates, the model exhibits the coexistence where the positive steady state is not necessarily unique \cite{ElhajjiIJB2018,ElhajjiJBD2009}.
Finally, the third type is the mutual inhibitory relationship which is characterized by the fact that each species inhibits the growth of the other species.
In \cite{ElhajjiIJB2018}, it was proved that the coexistence of two species is impossible which confirms the CEP.

Here, we focus attention on the third type of interaction which describes the mutual inhibitory relationship between species.
The present paper studies the competition model of two microbial species on a single nutriment \cref{ModelDEgaux} but with distinct removal rates.
By comparing with the existing literature, our study provides an extension of the results in \cite{ElhajjiIJB2018} by considering the mortality of each species.

In the first step, we want to understand the effect of decay on the behavior of the density-dependent model \cref{ModelDEgaux} considered in \cite{ElhajjiIJB2018}.
As is well known in the literature, the addition of mortality terms of the species can cause changes in the behaviors of the system.
For instance in predator-prey models in the chemostat, stable limit cycles \cite{BoerMB1998,MtarDcdsB2022} or multiple chaotic attractors \cite{KooiDCDISB2003} can be found when the mortality of species not neglected.
So, can the simultaneous presence of mortality and interspecific interference induce the stable coexistence of competitors in a chemostat?

In Fekih-Salem et al. \cite{FekihMB2017}, a model of competition of two species with inter and intraspecific density-dependent growth functions and distinct removal rates has been considered.
However, the particular case with only interspecific density-dependent growth rates has not been studied theoretically or numerically.

Here, we provide a complete mathematical and numerical analysis of the interspecific density-dependent model \cref{ModelDEgaux} by considering distinct removal rates for each species.
By allowing a general class of interspecific density-dependent growth functions, we describe the multiplicity of all steady states and their necessary and sufficient conditions of existence and local stability according to the operating parameters $S_{in}$ and $D$.

In the second step, our principal objective in this work is to provide an important tool which is the operating diagram of the system.
It provides a more global vision of the asymptotic behavior of solutions of the system according to the two operating parameters. 
On the other hand, the two-parameter diagram is also obtained by the numerical continuation method using MATCONT software \cite{MATCONT2023}. This method allows us to detect other types of bifurcations according to two parameters that cannot be detected theoretically.
Furthermore, we study the one-parameter bifurcation diagram with respect to the dilution rate $D$ showing the bifurcations of all steady states. 
Finally, through numerical simulations, we validate our theoretical results.
All these results have not been done in the literature.

The paper is structured as follows. 
First, we present in \Cref{Sec-SectMM} the mathematical model and some general hypotheses about the growth functions.
In \Cref{Sec-AnalMod}, we define some notations and functions used in this paper.
Then, we determine the existence, uniqueness, and local stability conditions of all steady states according to the dilution rate and the input concentration of the substrate.
Next, in \Cref{Sec-DO}, we describe theoretically the operating diagram of the system and then numerically using MATCONT.
The study of the one-parameter bifurcation with respect to the dilution rate and the numerical simulations are reported in \Cref{Sec-DB}.
Some conclusions are drawn in \Cref{Sec-Conc}.
\Cref{AppendixA} is devoted to the proofs of results presented in \Cref{Sec-AnalMod}.
All the values of parameters used throughout this paper are provided in \Cref{AppendixB}.
\section{Mathematical model}	                    \label{Sec-SectMM}
Taking into account the death rates of two species, the model is described by the following system of ordinary differential equations:
\begin{equation}                                \label{ModelYi}
\left\{\begin{array}{lll}
\dot{S}   &=& D(S_{in}-S) - \mu_1(S,X_2)X_1 - \mu_2(S,X_1)X_2,\\ 
\dot{X}_1 &=& (Y_1\mu_1(S,X_2)-D_1)X_1, \\
\dot{X}_2 &=& (Y_2\mu_2(S,X_1)-D_2)X_2,
\end{array}\right.
\end{equation}
where $Y_1 \leq 1 $ and $Y_2\leq 1$ are yield coefficients.
For $i=1,2$, $j=1,2$, $i \neq j$, $\mu_i$ is the density-dependent growth rate of the species $X_i$ that is assumed to be increasing in the variables $S$ and decreasing on $X_j$; 
$D_i$ is the removal rate of species $X_i$ and can be modeled as in \cite{FekihSIADS2021} by
$$
D_i=\alpha_i D+a_i, 
$$
where the coefficient $\alpha_i\in[0, 1]$ is a parameter for decoupling hydraulic retention time and solid retention time;
$a_i$ is the non-negative death rate of the species $X_i$.

To ease the mathematical analysis, we first can rescale system \cref{ModelYi} using the following change of variables $x_i=X_i/Y_i$. Thus, system \cref{ModelYi} becomes
\begin{equation}                              \label{DDInterSpecMod}
\left\{\begin{array}{lll}
\dot{S}   &=& D(S_{in}-S)- f_1(S,x_2)x_1-f_2(S,x_1)x_2,\\
\dot{x}_1 &=& (f_1(S,x_2)-D_1)x_1, \\ 
\dot{x}_2 &=& (f_2(S,x_1)-D_2)x_2,
\end{array}\right.
\end{equation}
where $f_1$ and $f_2$ are defined by 
$$
f_1(S,x_2)=Y_1\mu_1(S,Y_2x_2) \quad \mbox{and} \quad f_2(S,x_1)=Y_2\mu_2(S,Y_1x_1).
$$
In what follows, we consider model \cref{DDInterSpecMod}.
For $i=1,2$, $j=1,2$, $i\neq j$, we make the following general hypotheses on the growth rate $f_i: \mathbb{R}_+^2 \longrightarrow \mathbb{R}_+$ which is continuously differentiable.
\begin{hypothesis}                                        \label{hyp1}
No growth can occur for species $x_j$ without the substrate $S$ : 
$f_i(0,x_j)=0$, for all  $x_j \geq 0$.
\end{hypothesis}
\begin{hypothesis}                                        \label{hyp2}
The growth rate of species $x_i$ is increasing by the substrate $S$ but is inhibited by the species $x_j$ :
$\ds{\frac{\partial f_i}{\partial S}(S,x_j)>0}$ and $\ds{ \frac{\partial f_i}{\partial x_j}(S,x_j)<0}$ for all $S> 0$ and $x_j >0$.
\end{hypothesis}

As model \cref{DDInterSpecMod} represents a biological system, we require all solutions to be nonnegative and bounded for all time when given relevant initial conditions.
Thus, solutions of system \cref{DDInterSpecMod} verify the following result. Its proof uses classical tools, see for instance \cite{MtarDcdsB2022}.
\begin{proposition}
Assume that \cref{hyp1,hyp2} hold.
For any nonnegative initial condition, the solution of system \cref{DDInterSpecMod} exists for all $t \geq 0$, remains nonnegative and is bounded. Let $D_{\min}=\min (D,D_1,D_2)$. The set
$$
\Omega =\left\{ (S,x_1,x_2) \in \mathbb{R}_+^3 : S+x_1+x_2\leq D S_{in}/D_{\min} \right\}
$$
is positively invariant and is a global attractor for the dynamics \cref{DDInterSpecMod}.
\end{proposition}
\section{Existence and stability of steady states}	\label{Sec-AnalMod}
Our aim in this section is to determine the existence and local stability conditions of all steady states of model \cref{DDInterSpecMod}. For convenience, we use the abbreviation LES for Locally Exponentially Stable steady states. 
In all figures, the straight line $\delta$ is defined by the equation $D_1x_1/D+D_2x_2/D=S_{in}$ and the curves $\gamma_1$ and $\gamma_2$ are in blue and red, respectively.
In what follows, we need to define the set 
\begin{equation}                                   \label{SetM}
M := \left\{ (x_1,x_2) \in \mathbb{R}_+^2 : D_1x_1/D+D_2x_2/D \leq S_{in} \right\}
\end{equation}
and the following quantities:
\begin{equation}                                \label{ExprEFGH}
\ds{ E= \frac{\partial f_1}{\partial S},\quad F= \frac{\partial f_2}{\partial S},\quad G= -\frac{\partial f_1}{\partial x_2},\quad H= -\frac{\partial f_2}{\partial x_1}}.
\end{equation}
In \cref{ExprEFGH}, we have used the opposite sign of the partial derivatives $G=-\partial f_1/\partial x_2$ and $H=-\partial f_2/\partial x_1$, such that all constants involved in the computation become positive.

The steady states of \cref{DDInterSpecMod} are the solutions of the set of equations
\begin{equation}                         \label{SystSS}
\left\{\begin{array}{lll}
0 &=&  D(S_{in}-S)-f_1(S,x_2)x_1-f_2(S,x_1)x_2,\\
0 &=&  (f_1(S,x_2)-D_1)x_1, \\
0 &=&  (f_2(S,x_1)-D_2)x_2.
\end{array}\right.
\end{equation}
In the following, we will show that system \cref{DDInterSpecMod} can have at most four steady states labeled below:
\begin{itemize}
\item $\mathcal{E}_0=(S_{in},0,0)$: the washout of two species $(x_1=x_2=0)$.
\item $\mathcal{E}_1=(\tilde{S}_1,\tilde{x}_1,0)$: only the first species $x_1$ survives $(x_1>0,\: x_2=0)$.
\item $\mathcal{E}_2=(\tilde{S}_2,0,\tilde{x}_2)$: only the second species $x_2$ survives $(x_1=0,\: x_2>0)$.
\item $\mathcal{E}^*=(S^*,x_1^*,x_2^*)$: the coexistence of two species $(x_1>0, \:x_2>0)$.
\end{itemize}

\begin{figure}[!ht]
\setlength{\unitlength}{1.0cm}
\begin{center}
\begin{picture}(6.3,5)(0,0)
\put(-3.3,0){\rotatebox{0}{\includegraphics[width=4cm,height=5cm]{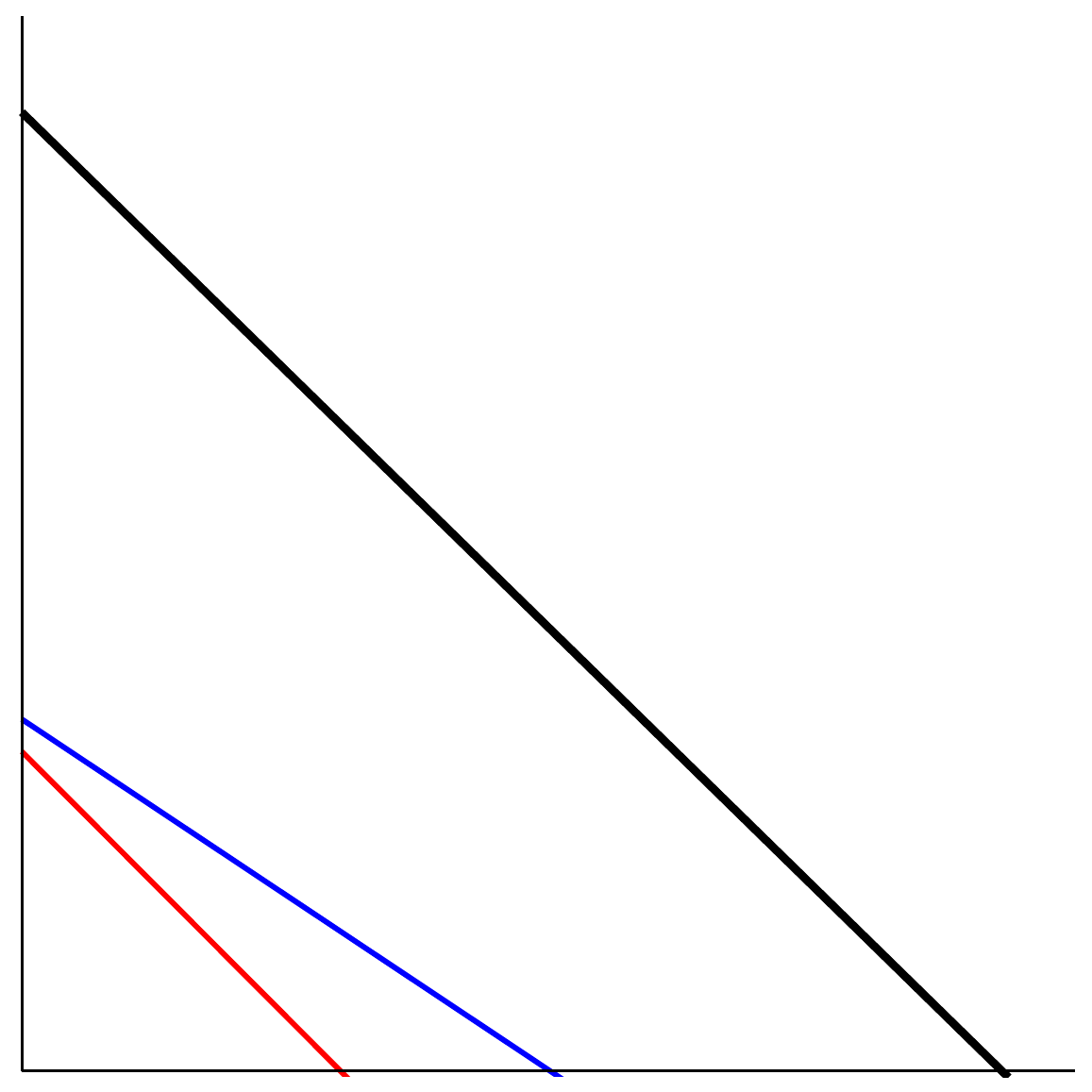}}}
\put(1.15,0){\rotatebox{0}{\includegraphics[width=4cm,height=5cm]{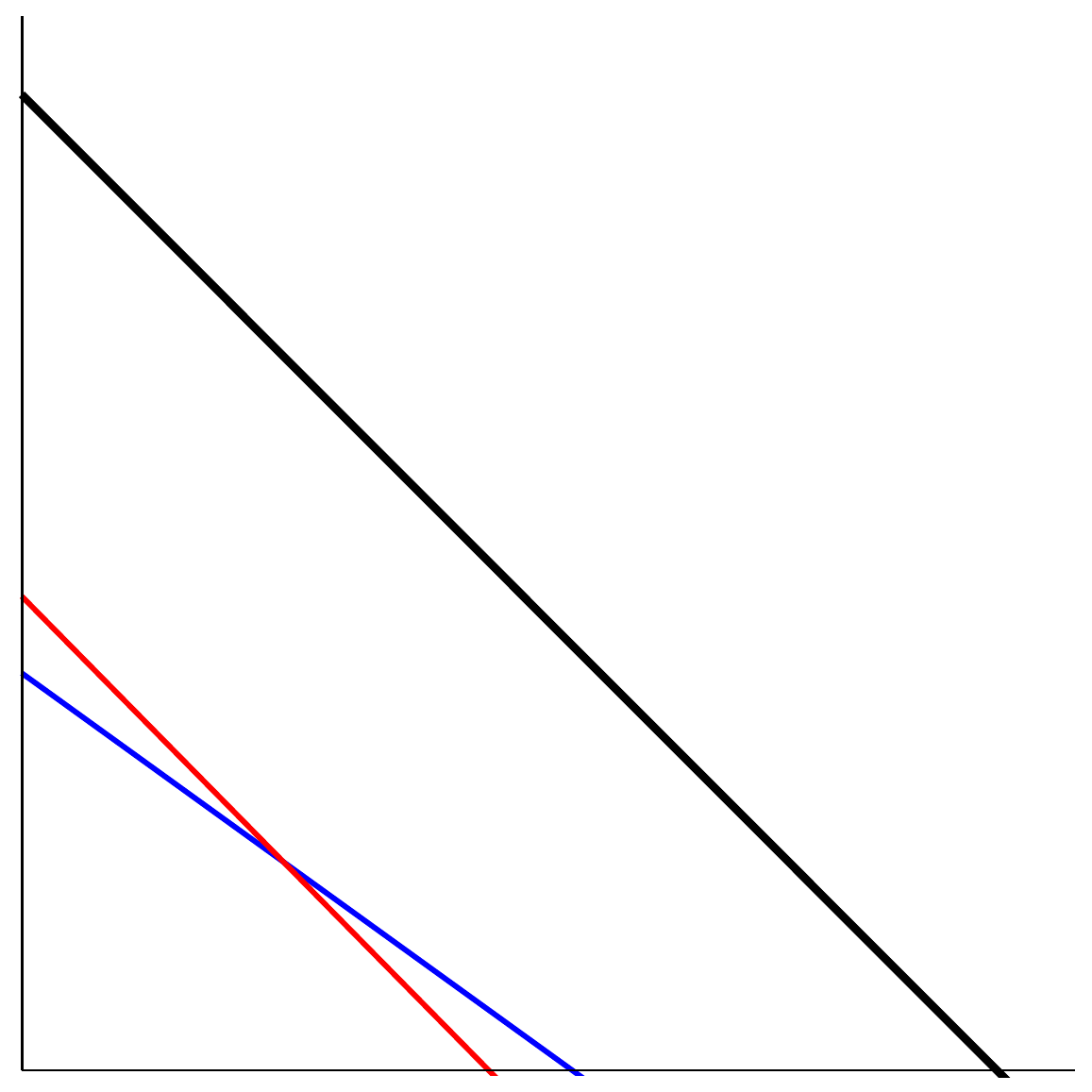}}}
\put(5.8,0){\rotatebox{0}{\includegraphics[width=4cm,height=5cm]{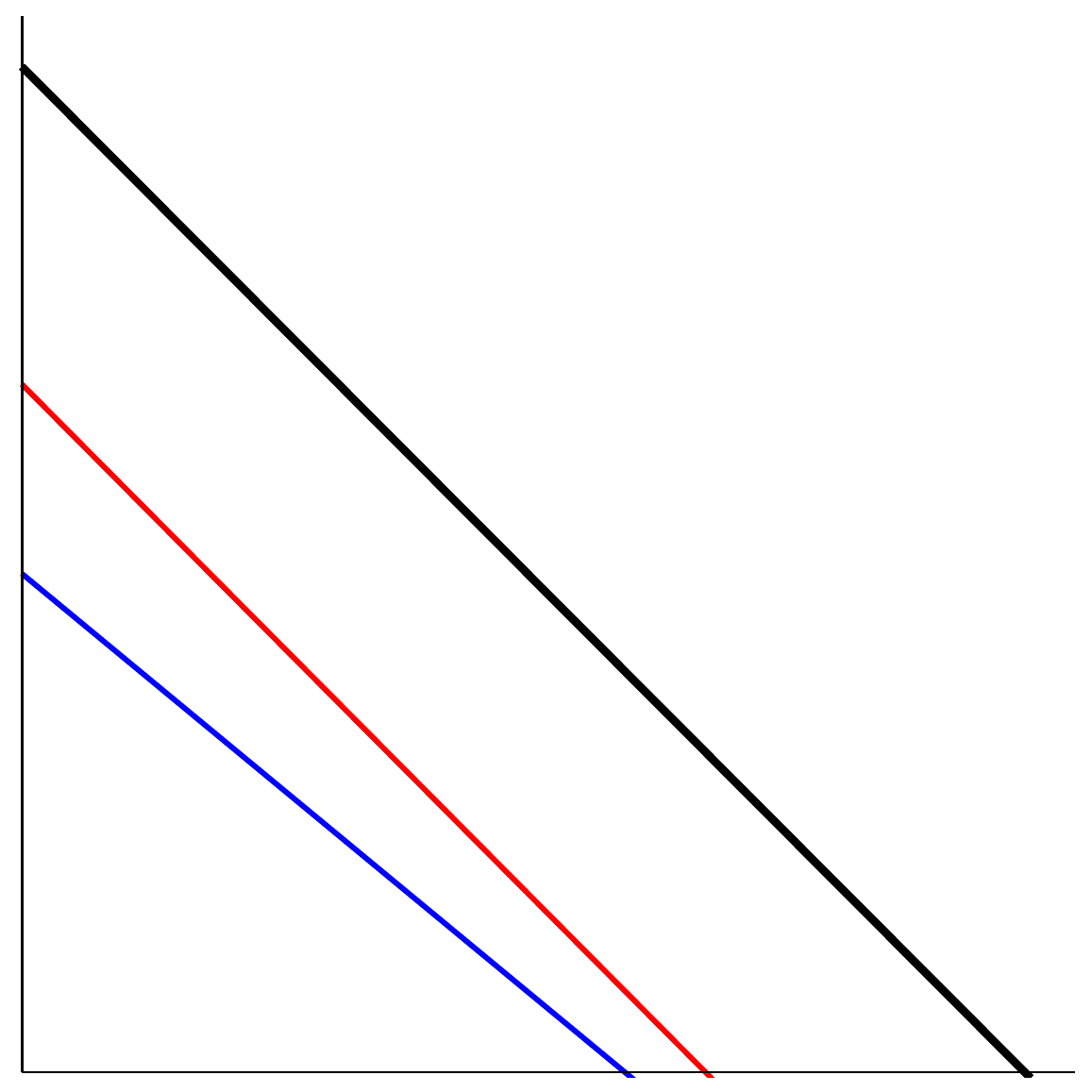}}}
\put(-1.8,4.8){{\sc (a)}}
\put(-3.15,4.65){{\sc  $x_2$}}
\put(-4.12,4.4){{\sc  $\frac{D}{D_2}S_{in}$}}
\put(-3.65,1.7){{\sc  $\bar{x}_2$}}
\put(-3.65,1.4){{\sc  $\tilde{x}_2$}}
\put(-2.3,1){{\sc {\color{blue}$\gamma_1$  }}}
\put(-2.35,0.5){{\sc {\color{red}$\gamma_2$  }}}
\put(-1.7,2.85){{\sc $\delta$  }}
\put(-2.2,-0.2){{\sc  $\bar{x}_1$}}
\put(-1.4,-0.2){{\sc  $\tilde{x}_1$}}
\put(0,-0.25){{\sc  $\frac{D}{D_1}S_{in}$}}
\put(0.4,0.2){{\sc  $x_1$}}
\put(-3.18,0.2){{\sc {\color{blue} $\mathcal{E}_0$}}}
\put(-3.27,0.02){{\sc {\color{blue} $\bullet$}}}
\put(-1.3,0.2){{\sc {\color{red} $\mathcal{E}_1$}}}
\put(-1.35,0.02){{\sc {\color{red} $\bullet$}}}
\put(-3.2,1.05){{\sc {\color{blue} $\mathcal{E}_2$}}}
\put(-3.28,1.45){{\sc {\color{blue} $\bullet$}}}
\put(2.8,4.8){{\sc (b)}}
\put(1.3,4.65){{\sc  $x_2$}}
\put(0.3,4.4){{\sc  $\frac{D}{D_2}S_{in}$}}
\put(2.9,0.5){{\sc {\color{blue}$\gamma_1$}}}
\put(1.75,1.7){{\sc {\color{red}$\gamma_2$  }}}
\put(2.8,2.7){{\sc $\delta$  }}
\put(3.2,-0.2){{\sc  $\tilde{x}_1$}}
\put(4.35,-0.25){{\sc  $\frac{D}{D_1}S_{in}$}}
\put(2.85,-0.2){{\sc  $\bar{x}_1$}}
\put(0.8,2.2){{\sc  $\tilde{x}_2$}}
\put(0.85,1.9){{\sc  $\bar{x}_2$}}
\put(4.85,0.2){{\sc  $x_1$}}
\put(1.27,0.2){{\sc {\color{blue} $\mathcal{E}_0$}}}
\put(1.17,0.02){{\sc {\color{blue} $\bullet$}}}
\put(3.25,0.2){{\sc {\color{red} $\mathcal{E}_1$}}}
\put(3.2,0.02){{\sc {\color{red} $\bullet$}}}
\put(1.3,2.3){{\sc {\color{red} $\mathcal{E}_2$}}}
\put(1.15,2.2){{\sc {\color{red} $\bullet$}}}
\put(2.3,1.15){{\sc {\color{blue} $\mathcal{E}^*$}}}
\put(2.1,1){{\sc {\color{blue} $\bullet$}}}
\put(7.6,4.8){{\sc (c)}}
\put(5.95,4.7){{\sc  $x_2$}}
\put(4.95,4.5){{\sc  $\frac{D}{D_2}S_{in}$}}
\put(5.48,3.2){{\sc  $\tilde{x}_2$}}
\put(5.48,2.3){{\sc  $\bar{x}_2$}}
\put(6.75,1.5){{\sc ${\color{blue}\gamma_1}$}}
\put(7.32,1.55){{\sc {\color{red}$\gamma_2$  }}}
\put(7.6,2.7){{\sc $\delta$  }}
\put(7.95,-0.2){{\sc  $\tilde{x}_1$}}
\put(9,-0.25){{\sc  $\frac{D}{D_1}S_{in}$}}
\put(8.36,-0.2){{\sc  $\bar{x}_1$}}
\put(9.55,0.2){{\sc  $x_1$}}
\put(5.9,0.2){{\sc {\color{blue} $\mathcal{E}_0$}}}
\put(5.8,0.02){{\sc {\color{blue} $\bullet$}}}
\put(7.52,0.25){{\sc {\color{blue} $\mathcal{E}_1$}}}
\put(7.99,0.02){{\sc {\color{blue} $\bullet$}}}
\put(5.95,3.35){{\sc {\color{red} $\mathcal{E}_2$}}}
\put(5.8,3.2){{\sc {\color{red} $\bullet$}}}
\end{picture}
\end{center}
\caption{Relative positions of $\tilde{x}_i$ and $\bar{x}_i$, $i=1,2$: (a) Case 1 of \cref{Cases123} no intersection, (b) Case 2 of \cref{Cases123} a unique intersection, (c) Case 3  of \cref{Cases123} no intersection. The red [resp. blue] color is for the LES [resp. unstable] steady state.}\label{Fig-CesesPos}
\end{figure}
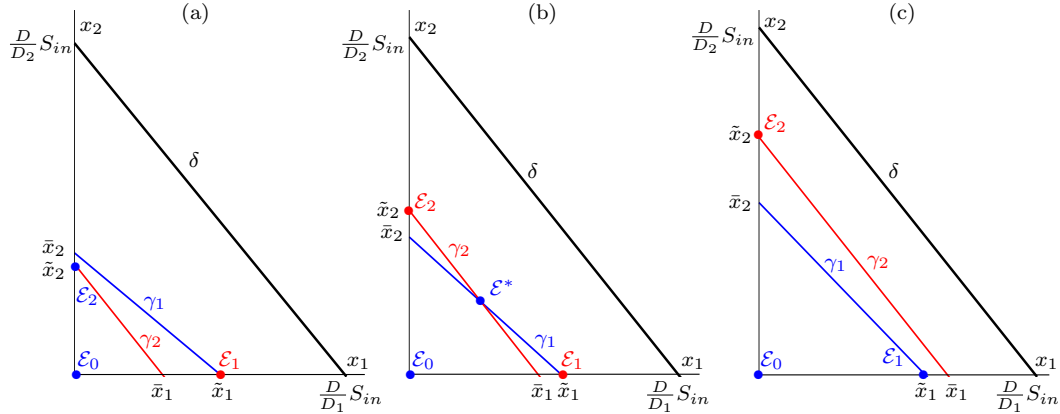

To determine the steady states of \cref{DDInterSpecMod} under \cref{hyp1,hyp2}, we need to define some notations and functions in \cref{Tab-BECAN} and \cref{Lem-Fi}, respectively. The proof of \cref{Lem-Fi} is similar to that of \cite{MtarCari2022}.
\begin{table}[ht]
\caption{Break-even concentrations and notations when $i=1,2$, $j=1,2$, $i\neq j$.} \label{Tab-BECAN}	
\begin{center}	                  
\begin{tabular}{ll}	
\begin{tabular}{l}
$\lambda_i(D)$  
\end{tabular}
&\begin{tabular}{l}
$S=\lambda_i(D)$ is the solution of $f_i(S,0)=\alpha_iD+a_i$   \\
It is defined for $D<(f_i(+\infty,0)-a_i)/\alpha_i$
\end{tabular}	
\\	  \hline 
\begin{tabular}{l}
$\tilde{x}_i$  
\end{tabular}
&\begin{tabular}{l}	
$\tilde{x}_i$ is the unique solution of $f_i(S_{in}-D_ix_i/D,0)=D_i$\\
$\tilde{x}_i=D(S_{in}-\lambda_i(D))/D_i$	
\end{tabular}
\\	  \hline  	
\begin{tabular}{l}
$\bar{x}_i$ \\
\end{tabular}
&\begin{tabular}{l}
$\bar{x}_i$ is the unique solution of $f_j(S_{in}-D_ix_i/D,x_i)=D_j$
\end{tabular}
\end{tabular}				
\end{center}
\end{table}
\begin{lemma}                             \label{Lem-Fi}
\begin{enumerate}[leftmargin=*] 
    \item Let $S_{in}>\lambda_1(D)$. The equation $f_1(S_{in}-D_1x_1/D-D_2x_2/D,x_2)=D_1$ defines a smooth decreasing function
$$
\begin{array}{ccccl}
F_1  & :  &  [0,\tilde{x}_1]  & \longrightarrow & [0,\bar{x}_2] \\
     &    &   x_1             & \longmapsto     &  F_1(x_1)=x_2,
\end{array}
$$
such that $F_1(\tilde{x}_1)=0$, $F_1(0)=\bar{x}_2$ where $\tilde{x}_1$ and $\bar{x}_2$ are defined in \cref{Tab-BECAN} with 
\begin{equation}                      \label{InegF1prim}
F'_1(x_1)=-\frac{D_1E}{D_2E+DG}<0, \quad \mbox{for all} \quad x_1 \in [0,\tilde{x}_1].
\end{equation}
Furthermore, the graph $\gamma_1$ of $F_1$ lies in the interior of $M$ for all $x_1 \in (0,\tilde{x}_1)$ (see \cref{Fig-CesesPos}).
\item Suppose that $S_{in}>\lambda_2(D)$.
The equation $f_2(S_{in}-D_1x_1/D-D_2x_2/D,x_1)=D_2$ defines a smooth decreasing function
$$
\begin{array}{ccccl}
F_2  & :  &  [0,\bar{x}_1]  & \longrightarrow & [0,\tilde{x}_2] \\
     &    &   x_1             & \longmapsto     &  F_2(x_1)=x_2,
\end{array}
$$
such that $F_2(\bar{x}_1)=0$, $F_2(0)=\tilde{x}_2$ where $\bar{x}_1$ and $\tilde{x}_2$ are defined in \cref{Tab-BECAN} with
\begin{equation}                      \label{InegF2prim}
F_2'(x_1)=-\frac{D_1F+DH}{D_2F}<0, \quad \mbox{for all} \quad x_1 \in [0,\bar{x}_1].
\end{equation}
In addition, the graph $\gamma_2$ of $F_2$ lies in the interior of $M$ for all $x_1 \in (0,\bar{x}_1)$ (see \cref{Fig-CesesPos}).
\end{enumerate}
\end{lemma}

Note that $\tilde{x}_i$ and $\bar{x}_i$, $i=1,2$ defined in \cref{Tab-BECAN} represent the coordinates of the intersections of the curves $\gamma_1$ and $\gamma_2$ with the coordinates axes.
According to their relative positions, there exist three cases that must be distinguished (see \cref{Fig-CesesPos}):
\begin{equation}                          \label{Cases123}
\begin{split}
\mbox{ Case ~1}:~  & \bar{x}_1<\tilde{x}_1 \quad \mbox{and} \quad \tilde{x}_2<\bar{x}_2.\\
\mbox{ Case ~2}:~  & \bar{x}_1<\tilde{x}_1 \quad \mbox{and} \quad \bar{x}_2<\tilde{x}_2.\\
\mbox{ Case ~3}:~  & \tilde{x}_1<\bar{x}_1 \quad \mbox{and} \quad \bar{x}_2<\tilde{x}_2.
\end{split}
\end{equation}
\begin{lemma}                                \label{Lemma-EquivXiFi}
When \cref{hyp1,hyp2} hold and $S_{in}>\lambda_i(D)$, for $i=1,2$, we have the following equivalences:
\begin{equation}                                   \label{EquivXiFi}
\begin{split}
\bar{x}_1<\tilde{x}_1 \quad \Leftrightarrow & \quad f_2(\lambda_1(D), D(S_{in}-\lambda_1(D))/D_1) < D_2,\\
\bar{x}_2<\tilde{x}_2 \quad \Leftrightarrow & \quad f_1(\lambda_2(D), D(S_{in}-\lambda_2(D))/D_2) < D_1.
\end{split}
\end{equation}
\end{lemma}  
The components of these steady states and the conditions of their existence according to the operating parameters are given by the following proposition:
\begin{proposition}                   \label{Prop-Exis}
Assume that \cref{hyp1,hyp2} hold.
The components of all steady states of system \cref{DDInterSpecMod} as well as their necessary and sufficient conditions of existence are given in \cref{Tab-ExisSS}. All steady states are unique, when they exist.
\begin{table}[ht]
\caption{Components of all steady states of model \cref{DDInterSpecMod} with their necessary and sufficient conditions of existence.} \label{Tab-ExisSS}
\begin{center}
\begin{tabular}{ @{\hspace{1mm}}l@{\hspace{3mm}} @{\hspace{3mm}}l@{\hspace{2mm}} @{\hspace{1mm}}l@{\hspace{1mm}} }
                     &       Components                                    & Existence condition \\ 
$\mathcal{E}_0$ &  $S=S_{in}$, $x_1=x_2=0$                                 & Always exists \\ \hline 
$\mathcal{E}_1$ &  $\tilde{S}_1=\lambda_1(D)$, $x_1=\tilde{x}_1$, $x_2=0$  & $S_{in}>\lambda_1(D)$\\ \hline 
$\mathcal{E}_2$ &  $\tilde{S}_2=\lambda_2(D) $, $x_1=0$, $x_2=\tilde{x}_2$ & $S_{in}>\lambda_2(D)$\\ \hline 
$\mathcal{E}^*$ &  \begin{tabular}{l}                                  
				   $S^*=S_{in}-D_1x_1^*/D-D_2x_2^*/D$,	\\
				   $(x_1^*,x_2^*)$ is the unique solution \\
                      of $x_2=F_1(x_1)=F_2(x_1)$
			\end{tabular}             & \begin{tabular}{l}    
                                         $f_1(\lambda_2(D), D(S_{in}-\lambda_2(D))/D_2)<D_1$ and \\
                                         $f_2(\lambda_1(D), D(S_{in}-\lambda_1(D))/D_1)<D_2$
                                         \end{tabular}     
\end{tabular}
\end{center}
\end{table}
\end{proposition}
As an immediate corollary of \cref{Lemma-EquivXiFi,Prop-Exis}, we have the following result establishing the existence and the uniqueness of the coexistence steady state $\mathcal{E}^*$.
\begin{corollary}                        \label{Coroll-UnicEet}
Assume that \cref{hyp1,hyp2} hold and $S_{in}>\lambda_i(D)$ for $i=1,2$.
A coexistence steady state $\mathcal{E}^*$ exists if and only if Case 2 of \cref{Cases123} holds. It is unique when it exists.
\end{corollary}
\begin{remark}                                   \label{Rq-PosF12}
If a positive intersection point $(x_1,x_2)$ of $\gamma_1$ and $\gamma_2$ exists, then the tangent of $\gamma_1$ at point $(x_1,x_2)$ is always above the tangent of $\gamma_2$ at point $(x_1,x_2)$.
\end{remark}
The stability conditions of all steady states of \cref{DDInterSpecMod} are provided by the following result:
\begin{proposition}                            \label{Prop-Stb}
Under \cref{hyp1,hyp2}, the necessary and sufficient conditions for the local asymptotic stability of the steady states $\mathcal{E}_0$, $\mathcal{E}_1$, $\mathcal{E}_2$ and $\mathcal{E}^*$ are given in \cref{Tab-StbSS}.
\end{proposition}
\begin{table}[ht]
\caption{Necessary and sufficient conditions of local stability of all steady states of \cref{DDInterSpecMod} according to the operating parameters $S_{in}$ and $D$. \label{Tab-StbSS}} 
\begin{center}
\begin{tabular}{ @{\hspace{1mm}}l@{\hspace{3mm}} @{\hspace{3mm}}l@{\hspace{1mm}} }
                &     Stability condition                              \\
$\mathcal{E}_0$ &  $S_{in}<\min (\lambda_1(D),\lambda_2(D))$            \\ \hline 
$\mathcal{E}_1$ &  $f_2(\lambda_1(D), D(S_{in}-\lambda_1(D))/D_1)<D_2$  \\ \hline 
$\mathcal{E}_2$ &  $f_1(\lambda_2(D), D(S_{in}-\lambda_2(D))/D_2)<D_1$  \\ \hline 
$\mathcal{E}^*$ &   unstable whenever it exists      
\end{tabular}
\end{center}
\end{table}
From \cref{Prop-Exis,Prop-Stb}, we obtain the following corollary.
\begin{corollary}
\begin{itemize}[leftmargin=*]
\item $\mathcal{E}_0$ is LES if and only if $\mathcal{E}_1$ and $\mathcal{E}_2$ do not exist.
\item A positive steady state $\mathcal{E}^*$ exists if and only if $\mathcal{E}_1$ and $\mathcal{E}_2$ are LES.
\end{itemize}
\end{corollary}
\section{Operating diagram}	                           \label{Sec-DO}
The operating diagram is a tool to describe how a system behaves when all biological parameters are fixed and the operating (control) parameters are varied, as they are the most easily manipulated parameters in a chemostat.
This operating diagram is very useful to understand the process from the mathematical and biological point of view.
It is often constructed in the mathematical literature (see for instance \cite{AbdellatifMBE2016,BarDCDSB2020,DaliYoucefMBE2020,DaoudMMNP2018,FekihSIADS2021,FekihMB2017,HarmandBook2017,MtarIJB2021,MtarDcdsB2022,SariNonLinDyn2021,SariMB2016,SariMB2017}) and the biological literature (see \cite{KhedimAMM2018,SbarciogBE2010,WadeJTB2016,XuJTB2011}).
In the existing literature, this operating diagram have been studied with three different methods.
For more details on these various methods, the reader is referred to \cite{MtarDcdsB2022} and the references therein.

The aim of this section is to study theoretically and numerically the operating diagram of \cref{DDInterSpecMod} under variation of the operating parameters which are the concentration of the substrate in the feed bottle $S_{in}$ and the dilution rate $D$. 
From \cref{Tab-ExisSS,Tab-StbSS}, we define in \cref{Tab-SetCurv} the set $\Upsilon=\left\{ \Upsilon_1,~ \Upsilon_2,~ \Upsilon_1^{c},~ \Upsilon_2^{c}\right\}$. This set represents the boundaries of the different regions in $(S_{in}, D)$-plane.
\begin{table}[ht]
\caption{The set $\Upsilon$ of boundaries of the regions in the operating diagram and their corresponding colors used in \cref{FigDOMaple,Fig-DOMatc}. \label{Tab-SetCurv}} 
\begin{center}
\begin{tabular}{ @{\hspace{1mm}}l@{\hspace{3mm}} @{\hspace{3mm}}l@{\hspace{1mm}}  }
$\Upsilon_1=\left\{ (S_{in},D): S_{in} = \lambda_1(D) \right\}$              &      Black       \\ \hline 
$\Upsilon_2=\left\{ (S_{in},D): S_{in} = \lambda_2(D) \right\}$              &      Blue        \\ \hline 
$\Upsilon_1^{c}=\left\{(S_{in},D): \tilde{x}_1(S_{in},D)=\bar{x}_1(S_{in},D) \right\}$  &      Red     \\ \hline 
$\Upsilon_2^{c}=\left\{(S_{in},D): \tilde{x}_2(S_{in},D)=\bar{x}_2(S_{in},D) \right\}$  &    Magenta      \\ 
\end{tabular}
\end{center}
\end{table}
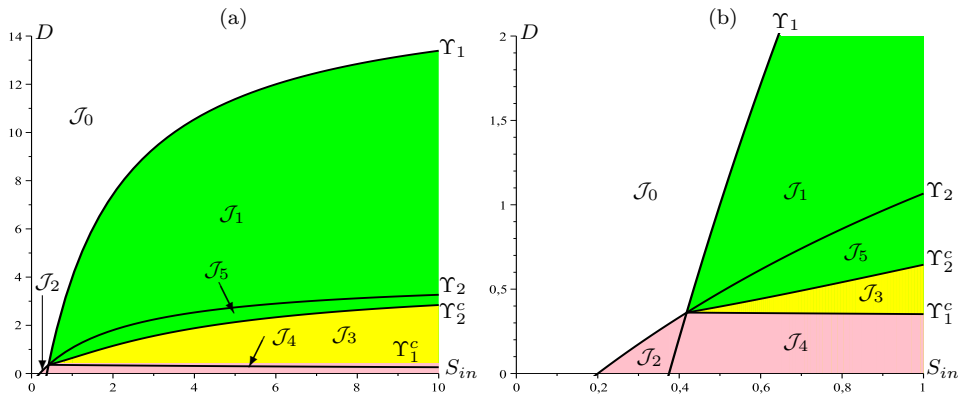
\begin{figure}[!ht]
\setlength{\unitlength}{1.0cm}
\begin{center}
\begin{picture}(6.3,5.1)(0,0)
\put(-3.4,0){\rotatebox{0}{\includegraphics[width=6cm,height=5cm]{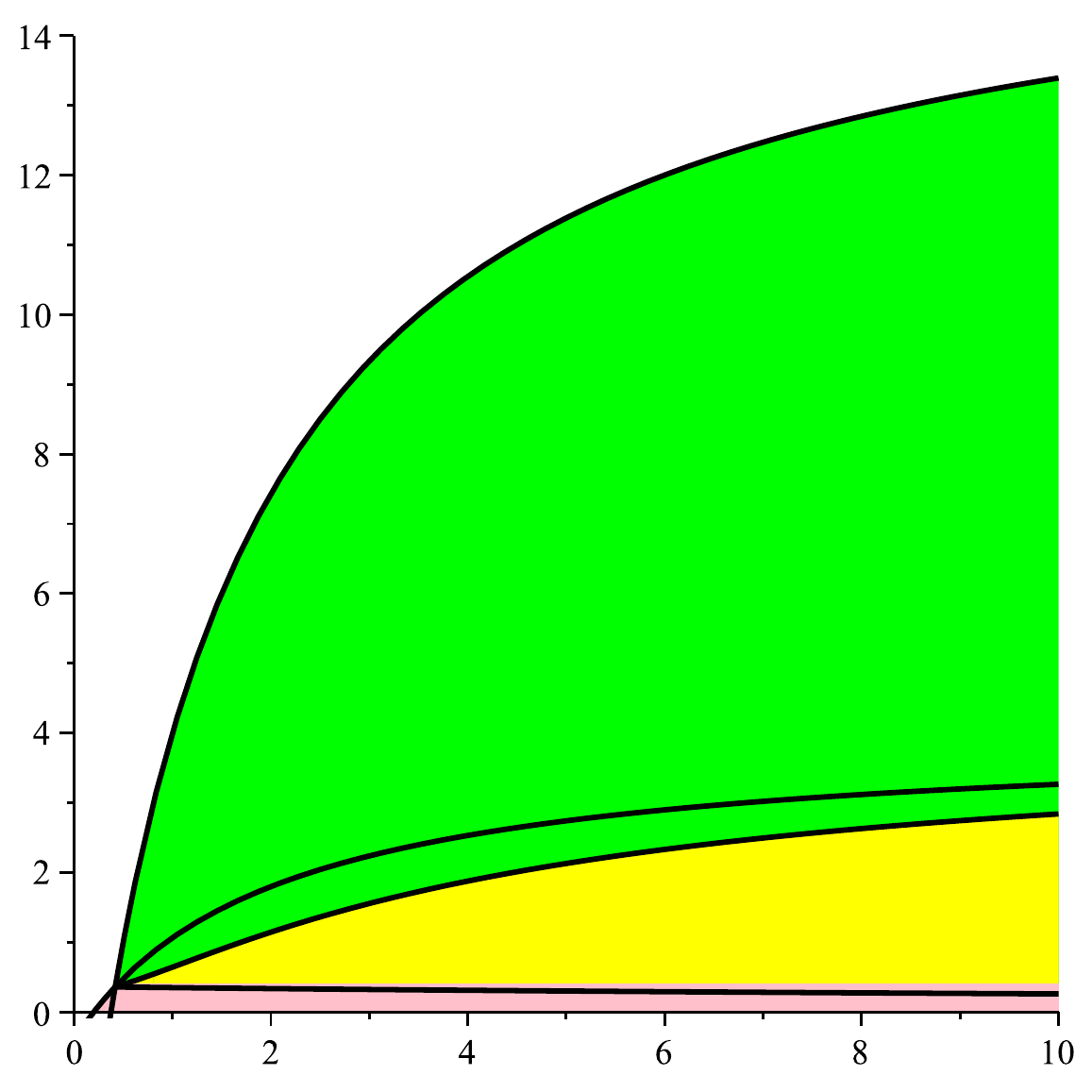}}}
\put(3,0){\rotatebox{0}{\includegraphics[width=6cm,height=5cm]{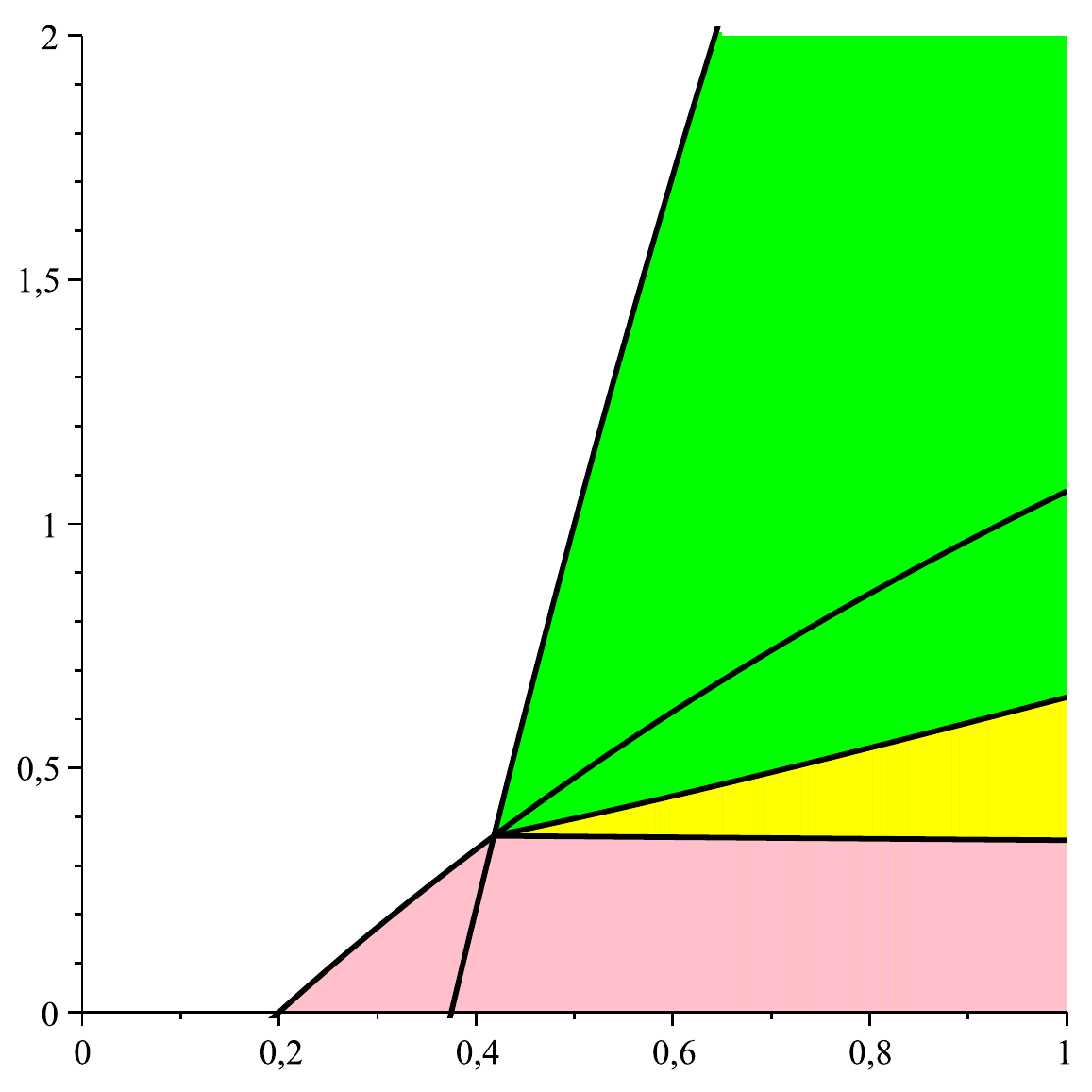}}}
\put(-0.5,5){\sc (a)}	
\put(2.5,0.35){{\sc  $S_{in}$}}
\put(-2.95,4.8){{\sc  $D$}}
\put(-0.5,2.4){{\sc  $\mathcal{J}_1$ }}
\put(-2.5,3.7){{\sc $\mathcal{J}_0$}}
\put(-0.7,1.65){{\sc  $\mathcal{J}_5$ }}
\put(-0.5,1.57){\sc \vector(1,-2){0.2}}
\put(1,0.8){{\sc  $\mathcal{J}_3$ }}
\put(0.2,0.77){{\sc  $\mathcal{J}_4$ }}
\put(0.1,0.85){\sc \vector(-1,-2){0.2}}
\put(-2.95,1.5){{\sc  $\mathcal{J}_2$}}
\put(-2.85,1.4){\sc \vector(0,-1){1}}
\put(2.42,4.6){{\sc  $\Upsilon_1$}}
\put(2.42,1.45){{\sc $\Upsilon_2$}}
\put(2.42,1.1){{\sc  $\Upsilon_2^{c}$}}
\put(1.8,0.6){{\sc  $\Upsilon_1^{c}$}}
\put(6,5){{\sc (b)}}
\put(8.9,0.35){{\sc  $S_{in}$}}
\put(3.5,4.8){{\sc  $D$}}
\put(8,1.35){{\sc  $\mathcal{J}_3$}}
\put(7.8,1.9){{\sc $\mathcal{J}_5$}}
\put(7,0.7){{\sc $\mathcal{J}_4$}}
\put(5,2.7){{\sc  $\mathcal{J}_0$}}
\put(7,2.7){{\sc $\mathcal{J}_1$}}
\put(5.05,0.5){{\sc  $\mathcal{J}_2$}}
\put(6.85,4.95){{\sc  $\Upsilon_1$}}
\put(8.9,2.7){{\sc  $\Upsilon_2$}}
\put(8.9,1.8){{\sc  $\Upsilon_2^{c}$}}
\put(8.9,1.1){{\sc $\Upsilon_1^{c}$}}
\end{picture}
\end{center}
\vspace{-0.3cm}
\caption{MAPLE: (a) operating diagram of \cref{DDInterSpecMod} in the $(S_{in},D)$-plane. (b) Magnification on the regions $\mathcal{J}_2$ and $\mathcal{J}_4$ when $(S_{in},D) \in [0,1]\times[0,2]$.} \label{FigDOMaple}
\end{figure}

In order to construct the operating diagram, we consider the following specific growth rates satisfying \cref{hyp1,hyp2}:
\begin{equation}                                     \label{SpeciFunc}
f_1(S,x_2)= \frac{m_1 S}{K_1+S+\beta_1x_2}, \quad f_2(S,x_1)= \frac{m_2S}{K_2+S+\beta_2x_1},
\end{equation}
where $m_1$ and $m_2$ denote the maximum growth rates;
$K_1$ and $K_2$ denote the Michaelis-Menten constants;
$\beta_1$ [resp. $\beta_2$] represents the inhibition factor due to $x_2$ [resp. $x_1$] for the growth of the species $x_1$ [resp. $x_2$].
Except for the operating parameters $S_{in}$ and $D$, all these parameters cannot be easily manipulated by the biologist.
Thus, we fix it as provided in \cref{Tab-Allpar}. 
Note that the construction of the operating diagram is similar for any other specific growth rate satisfying \cref{hyp1,hyp2}.

From \cref{Tab-ExisSS,Tab-StbSS} establishing the existence and local stability conditions of all steady states, we can state the following result determining theoretically the operating diagram where the various functions and the corresponding curves are defined in \cref{Tab-BECAN,Tab-SetCurv}.
\begin{proposition}
For the specific growth rates $f_1$ and $f_2$ defined in \cref{SpeciFunc} and the biological parameter values provided in \cref{Tab-Allpar},
the existence and the local stability of all steady states of system \cref{DDInterSpecMod} in the six regions $\mathcal{J}_k$, $k=0,\ldots,5$ of the operating diagram shown in \cref{FigDOMaple} are described in \cref{{Tab-DefReg}}.
\end{proposition}
\begin{table}[ht]
\caption{Definitions of the regions of the operating diagram in \cref{FigDOMaple,Fig-DOMatc} when $i=1,2$. The letter S [resp. U)] means stable [resp. unstable] steady state. The absence of a letter means that the corresponding steady state does not exist.\label{Tab-DefReg}} 
\begin{center}
\begin{tabular}{ @{\hspace{1mm}}l@{\hspace{1mm}} @{\hspace{1mm}}l@{\hspace{1mm}} @{\hspace{1mm}}l@{\hspace{1mm}} @{\hspace{1mm}}l@{\hspace{1mm}}  @{\hspace{1mm}}l@{\hspace{1mm}} @{\hspace{1mm}}l@{\hspace{1mm}} @{\hspace{1mm}}l@{\hspace{1mm}} @{\hspace{1mm}}l@{\hspace{1mm}} }
Region                        &      Condition $1$          & Condition $2$                               & Color  & $\mathcal{E}_0$ & $\mathcal{E}_1$ & $\mathcal{E}_2$ & $\mathcal{E}^*$  \\ 
$\mathcal{J}_0$ & $S_{in}>\min (\lambda_1(D),\lambda_2(D))$ &                                             & White  &    S            &                 &                 &                   \\ \hline                                              
$\mathcal{J}_1$ & $\lambda_1(D)<S_{in}<\lambda_2(D)$        & $\bar{x}_1(S_{in},D)<\tilde{x}_1(S_{in},D)$ & Green  &    U            &       S         &                 &                   \\ \hline  
$\mathcal{J}_2$ & $\lambda_2(D)<S_{in}<\lambda_1(D)$        & $\bar{x}_2(S_{in},D)<\tilde{x}_2(S_{in},D)$ & Pink   &    U            &                 &       S         &                   \\ \hline  
$\mathcal{J}_3$ & $S_{in}>\max(\lambda_1(D),\lambda_2(D))$  & $\bar{x}_i(S_{in},D)<\tilde{x}_i(S_{in},D)$ & Yellow &    U            &       S         &       S         &        U          \\ \hline  
$\mathcal{J}_4$ & $S_{in}>\max(\lambda_1(D),\lambda_2(D)) $ &  
\begin{tabular}{l}    
 $\bar{x}_2(S_{in},D)<\tilde{x}_2(S_{in},D)$, \\
$\bar{x}_1(S_{in},D)>\tilde{x}_1(S_{in},D)$
\end{tabular} 
                                                                                                          &  Pink  &    U            &       U         &       S         &                   \\ \hline  
$\mathcal{J}_5$  & $S_{in}>\max(\lambda_1(D),\lambda_2(D)) $& 
\begin{tabular}{l}    
$\bar{x}_2(S_{in},D)>\tilde{x}_2(S_{in},D)$, \\
$\bar{x}_1(S_{in},D)<\tilde{x}_1(S_{in},D)$
\end{tabular} 
                                                                                                          &  Green  &    U            &       S         &       U         &                   \\ \hline  
\end{tabular}
\end{center}
\end{table}

As shown in \cref{FigDOMaple}, the $(S_{in},D)$-plane is divided into six regions. 
The region $\mathcal{J}_0$ corresponds to the washout of two species.
The regions $\mathcal{J}_1$ and $\mathcal{J}_5$ [resp. $\mathcal{J}_2$ and $\mathcal{J}_4$] correspond to the competitive exclusion of the second [resp.  first] species. 
The region $\mathcal{J}_3$ corresponds to the bi-stability with convergence to one of the two boundary steady states.

In the following, we used MATCONT \cite{MATCONT2023} to numerically analyze the two-parameter diagram of system \cref{DDInterSpecMod} for the set of the biological parameter values of \cref{FigDOMaple}, which is provided in \cref{Tab-Allpar}.
Indeed, MATCONT \cite{MATCONT2023} is a MATLAB numerical continuation package for the bifurcation study of continuous and discrete parameterized systems of ODEs.
Thanks to the numerous features and functions, it allows us to compute curves of steady states and limit cycles as well as their local asymptotic behaviors and bifurcations as Branch Points (BP) or transcritical bifurcations, Limit Points (LP) or saddle-node (or fold) bifurcations, Hopf points (H), Limit Point of Cycles (LPC), Cusp bifurcations (CP), Bogdanov-Takens bifurcations, etc. 
All these curves are determined under variation of one or more system parameters by numerical continuation.
For more details, the reader is referred to \cite{DhoogeMCMD2008} and the references therein.
\begin{figure}[!ht]
\setlength{\unitlength}{1.0cm}
\begin{center}
\begin{picture}(6.3,10.9)(0,0)
\put(-4.9,0){\rotatebox{0}{\includegraphics[width=9cm,height=16.5cm]{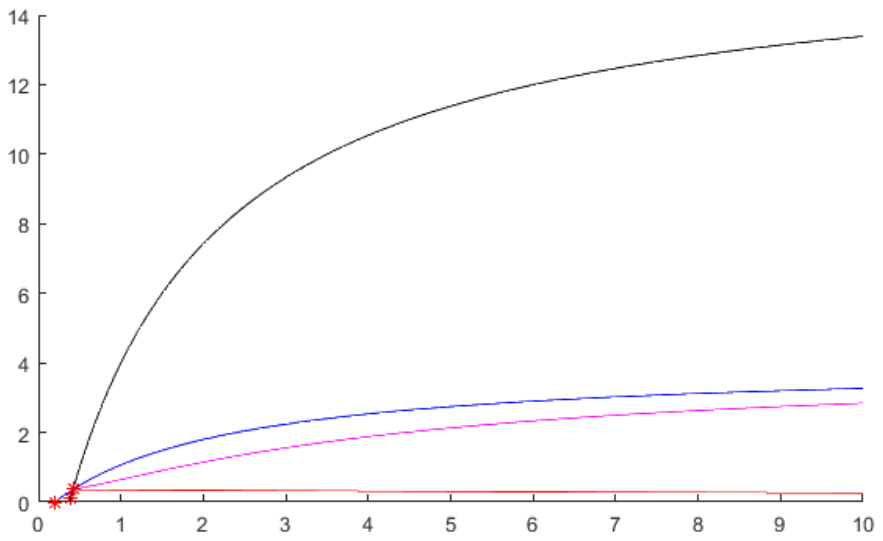}}}
\put(2.15,0){\rotatebox{0}{\includegraphics[width=8.8cm,height=16.5cm]{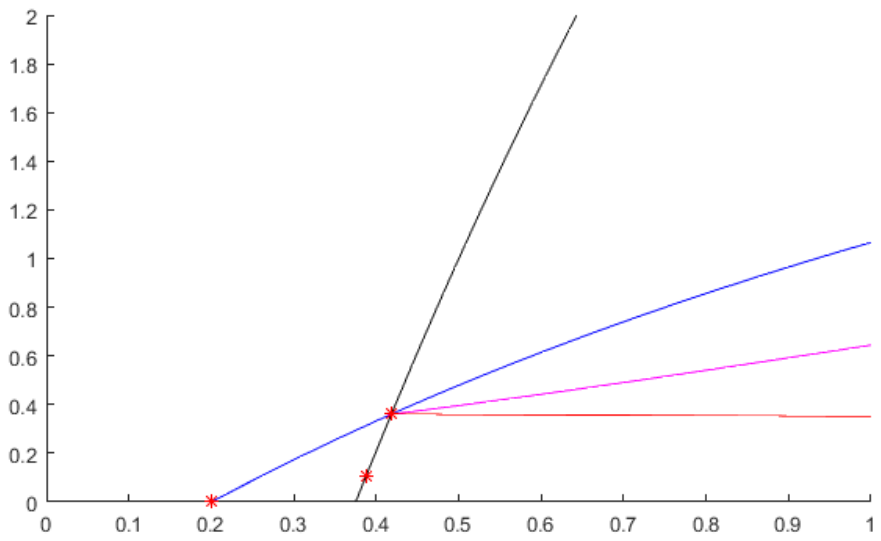}}}
\put(-0.5,10.6){{\sc (a)}}
\put(2.75,6.05){{\sc  $S_{in}$}}
\put(-3.2,10.6){{\sc  $D$}}
\put(-0.2,8.5){{\sc  $\mathcal{J}_1$ }}
\put(-2.5,10){{\sc $\mathcal{J}_0$}}
\put(0.85,7.5){{\sc  $\mathcal{J}_5$ }}
\put(1.1,7.4){\sc \vector(1,-2){0.2}}
\put(1,6.5){{\sc  $\mathcal{J}_3$ }}
\put(-0.1,6.5){{\sc  $\mathcal{J}_4$ }}
\put(-0.15,6.5){\sc \vector(-1,-2){0.2}}
\put(-3.13,7.2){{\sc  $\mathcal{J}_2$}}
\put(-3.08,7.1){\sc \vector(0,-1){1}}
\put(-2.55,7.1){{\sc  \texttt{ZH}}}
\put(-2.58,7.05){\sc \vector(-1,-2){0.4}}
\put(-3.15,5.75){{\sc  \texttt{BT}}}
\put(6.6,10.6){{\sc (b)}}
\put(9.65,6.05){{\sc  $S_{in}$}}
\put(3.8,10.6){{\sc  $D$}}
\put(9,7.05){{\sc  $\mathcal{J}_3$}}
\put(9,7.8){{\sc $\mathcal{J}_5$}}
\put(7.8,6.35){{\sc $\mathcal{J}_4$}}
\put(5,9){{\sc  $\mathcal{J}_0$}}
\put(7.8,9){{\sc $\mathcal{J}_1$}}
\put(5,6.9){{\sc  $\mathcal{J}_2$}}
\put(5.4,6.8){\sc \vector(1,-2){0.2}}
\put(4.7,6.2){{\sc  \texttt{BT}}}
\put(5.8,6.9){{\sc  \texttt{ZH}}}
\put(6.15,6.2){{\sc  \texttt{ZH}}}
\end{picture}
\end{center}
\vspace{-6cm}
\caption{MATCONT: (a) Operating diagram of \cref{DDInterSpecMod}. (b) Magnification of the two-parameter bifurcations at points ZH and BT when $(S_{in},D) \in [0,1]\times[0,2]$.}  \label{Fig-DOMatc}
\end{figure}

\cref{Fig-DOMatc} illustrates the operating diagram obtained numerically using MATCONT which is identical to one determined theoretically in \cref{FigDOMaple}. 
However, MATCONT can detect other two-parameter bifurcations as Zero-Hopf Point (ZH) and Bogdanov-Takens Point (BT).
These types of bifurcation are not detected in the theoretical study of the operating diagram.
In \cref{Fig-DOMatc}, the intersection point of all curves in the set $\Upsilon$ corresponds to a two-parameter bifurcation of type Zero-Hopf Point (ZH).
Moreover, MATCONT was able to detect a second Zero-Hopf Point (ZH) and a new two-parameter bifurcation of type Bogdanov-Takens (BT) with normal form coefficient $(a,b)=(2.56\: 10^{-6}, 0.0909)$.
\cref{Tab-VPDOMatc} summarizes the critical operating parameters and the state at BT and ZH bifurcations.
\begin{table}[ht]
\caption{Operating parameters and state for BT and ZH at the bifurcation points in \cref{Fig-DOMatc}. The abbreviation BT [resp. ZH] means a Bogdanov-Takens [resp. Zero-Hopf] bifurcation. \label{Tab-VPDOMatc}} 
\begin{center}
\begin{tabular}{ @{\hspace{1mm}}l@{\hspace{2mm}} @{\hspace{2mm}}l@{\hspace{2mm}} @{\hspace{2mm}}l@{\hspace{1mm}} }
Bifurcation & Parameter $(S_{in},D)$ & State $(S,x_1,x_2)$      \\  
ZH          & $(0.418289,0.361063)$  & $(0.418289,0,0)$         \\ \hline
BT          & $(0.2,2.9\:10^{-6})$   & $(0.2,0,0)$              \\ \hline  
ZH          & $(0.387348,0.10468)$   & $(0.387348,0,0)$         \\
\end{tabular}
\end{center}
\end{table}
\section{Bifurcation analysis and numerical simulations}	\label{Sec-DB}
In this section, we will analyze the one-parameter bifurcation diagram of system \cref{DDInterSpecMod}.
This diagram shows the various types of bifurcations and the asymptotic behavior of the process by varying the dilution rate $D$ as the bifurcating parameter where the input concentration $S_{in}$ is fixed.
In a similar way, we can study the one-parameter bifurcation diagram where $S_{in}$ is the bifurcation parameter.
Note that we have chosen $S_{in}=1$ to maximize the passage number through the various regions of the operating diagram in \cref{FigDOMaple,Fig-DOMatc}. 
The following result establishes the nature of all bifurcations that occur by passing from one region to another in the operating diagram in \cref{FigDOMaple,Fig-DOMatc}. 
\begin{proposition}       
For the specific growth rates $f_1$ and $f_2$ \cref{SpeciFunc} and the set of the biological parameter values provided in \cref{Tab-Allpar}, the bifurcation analysis of all steady states of system \cref{DDInterSpecMod} 
by crossing the various curves in the set $\Upsilon$ defined in \cref{Tab-SetCurv}, is summarized in \cref{Tab-TransBif}.
\begin{table}[ht]
\caption{Nature of bifurcations of steady states of \cref{DDInterSpecMod} by crossing to the surfaces of $\Upsilon$. The letter TB means a transcritical bifurcation. \label{Tab-TransBif}} 
\begin{center}
\begin{tabular}{ @{\hspace{1mm}}l@{\hspace{2mm}} @{\hspace{2mm}}l@{\hspace{2mm}}  @{\hspace{2mm}}l@{\hspace{2mm}}  @{\hspace{2mm}}l@{\hspace{1mm}} }
Transition &    $\Upsilon$       &  Bifurcation  & Steady state \\ 
$\mathcal{J}_0$ to $\mathcal{J}_1$ or $\mathcal{J}_2$ to $\mathcal{J}_4$ &$\Upsilon_1$   &  TB & $\mathcal{E}_0=\mathcal{E}_1$   \\ \hline 
$\mathcal{J}_0$ to $\mathcal{J}_2$ or $\mathcal{J}_1$ to $\mathcal{J}_5$ &$\Upsilon_2$   &  TB & $\mathcal{E}_0=\mathcal{E}_2$       \\ \hline 
$\mathcal{J}_3$ to $\mathcal{J}_4$ &$\Upsilon_1^{c}$ &  TB & $\mathcal{E}_1=\mathcal{E}^*$    \\ \hline 
$\mathcal{J}_5$ to $\mathcal{J}_3$ &$\Upsilon_2^{c}$ &  TB & $\mathcal{E}_2=\mathcal{E}^*$       
\end{tabular}
\end{center}
\end{table}
\end{proposition}

\begin{figure}[!ht]
\setlength{\unitlength}{1.0cm}
\begin{center}
\begin{picture}(6.3,5.3)(0,0)
\put(-3.15,0){\rotatebox{0}{\includegraphics[width=6cm,height=5cm]{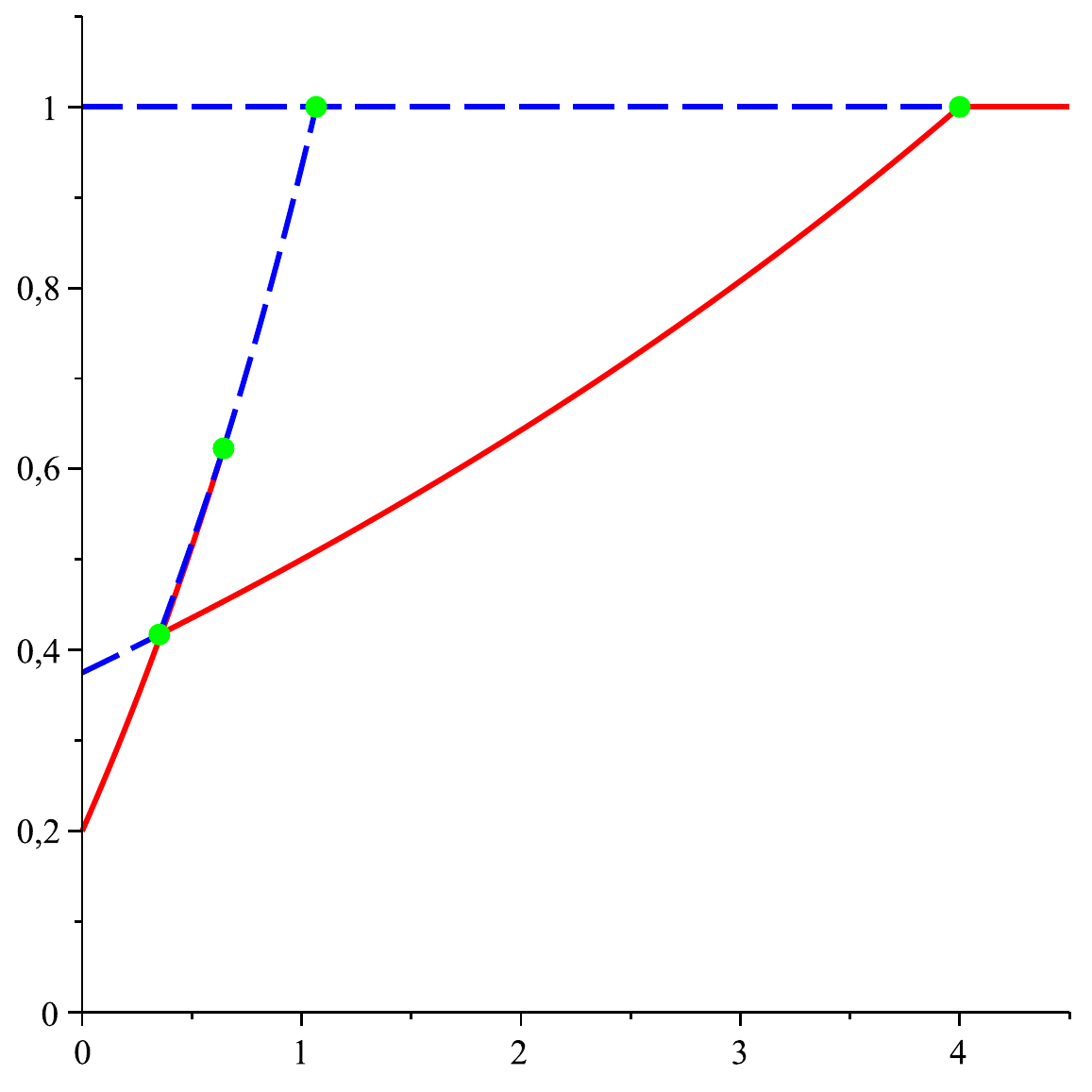}}}
\put(3.1,0){\rotatebox{0}{\includegraphics[width=6cm,height=5cm]{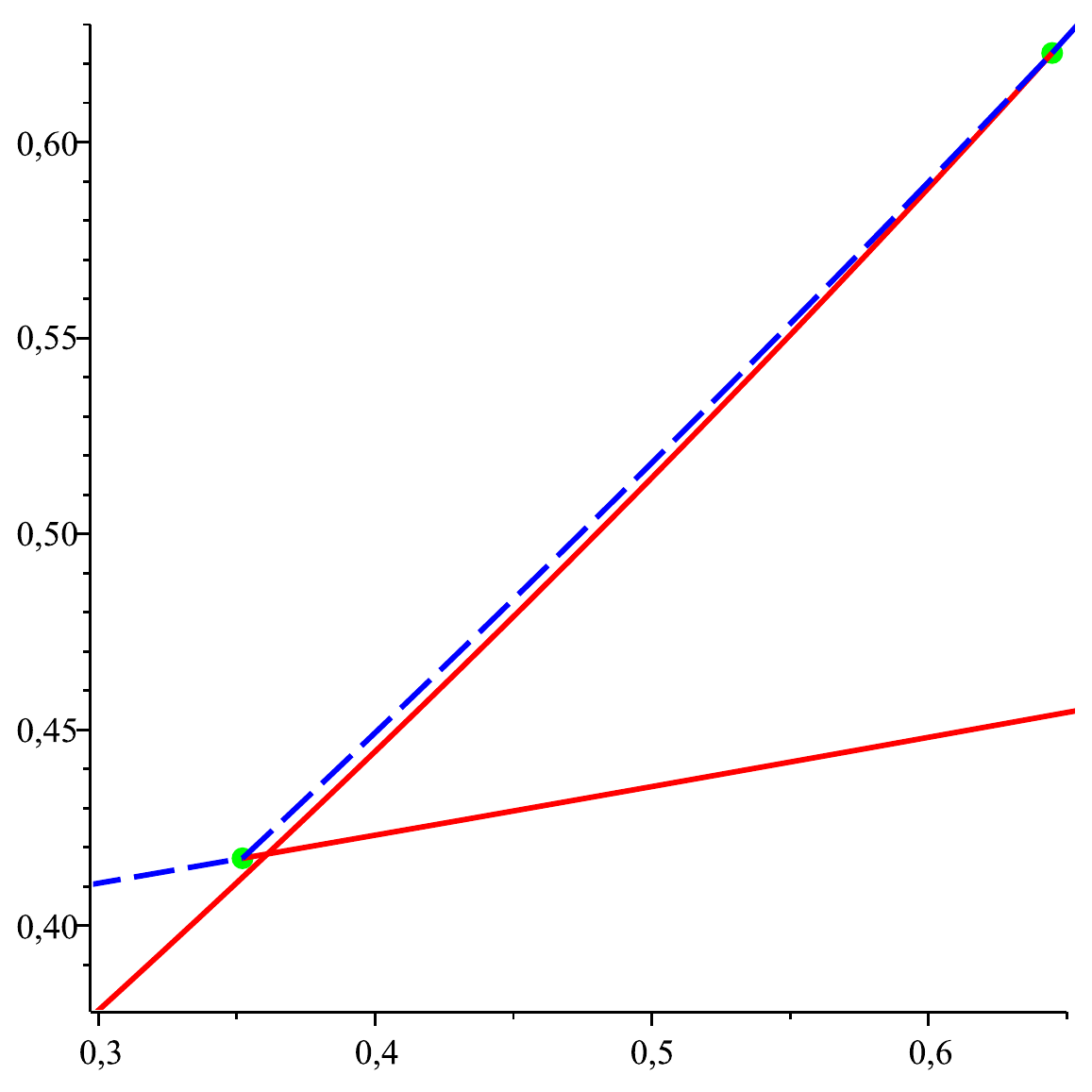}}}
\put(-0.3,4.9){\sc (a)}	
\put(2.75,0.36){{\sc  $D$}}
\put(-2.65,4.85){{\sc  $S$}}
\put(-1.2,4.6){{\sc {\color{blue}$\mathcal{E}_0$  }}}
\put(2.3,4.6){{\sc {\color{red}$\mathcal{E}_0$  }}}
\put(-0.5,3.15){{\sc  {\color{red}$\mathcal{E}_1$  } }}
\put(-1.5,3.85){{\sc {\color{blue}$\mathcal{E}_2$  }}}
\put(-2,2.5){{\sc  {\color{red}$\mathcal{E}_2$  }}}
\put(-2.45,1.45){{\sc  {\color{red}$\mathcal{E}_2$  }}}
\put(-2.65,2.12){{\sc  {\color{blue}$\mathcal{E}_1$  }}}
\put(-2.4,2.56){{\sc  {\color{blue}$\mathcal{E}^*$  }}}
\put(2.05,-0.05){{\sc  $\sigma_1$}}
\put(2.05,0.25){{\sc  $\dagger$}}
\put(-1.55,-0.05){{\sc  $\sigma_2$}}
\put(-1.55,0.25){{\sc  $\dagger$}}
\put(-2.05,-0.05){{\sc  $\sigma_3$}}
\put(-2.05,0.25){{\sc  $\dagger$}}
\put(-2.4,-0.05){{\sc  $\sigma_4$}}
\put(-2.4,0.25){{\sc  $\dagger$}}
\put(6,4.9){{\sc (b)}}
\put(9,0.38){{\sc  $D$}}
\put(3.65,4.85){{\sc  $S$}}
\put(6.9,2.7){{\sc {\color{red}$\mathcal{E}_2$  }}}
\put(6.7,3.15){{\sc {\color{blue}$\mathcal{E}^*$  }}}
\put(8.75,-0.05){{\sc  $\sigma_3$}}
\put(8.75,0.25){{\sc  $\dagger$}}
\put(4.3,-0.05){{\sc  $\sigma_4$}}
\put(4.3,0.25){{\sc  $\dagger$}}
\put(7.4,1.6){{\sc  {\color{red}$\mathcal{E}_1$  }}}
\put(3.75,1.1){{\sc  {\color{blue}$\mathcal{E}_1$  }}}
\end{picture}
\end{center}
\vspace{-0.2cm}
\caption{(a) One-parameter bifurcation diagram of \cref{DDInterSpecMod} in the variable $S$ with $D$ as the bifurcation parameter and $S_{in}=1$. (b) Magnification when $D \in [\sigma_3,\sigma_4]$. Blue dashed curves correspond to unstable steady states and red solid curves to LES steady states. The green solid diamonds represent the transcritical bifurcations.} \label{Fig-DB}
\end{figure}
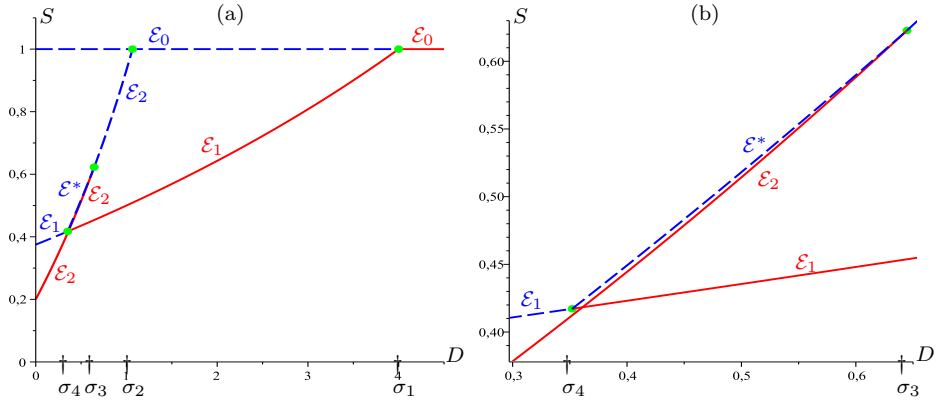

From \cref{Tab-TransBif} and as we shall see later, the passage from one region to another through the curves of the set $\Upsilon$ always corresponds to a transcritical bifurcation where steady states coalesce and change stability.
\Cref{Fig-DB} represents the one-parameter bifurcation diagram of system \cref{DDInterSpecMod} where the $\omega$-limit set is projected to the $S$ coordinate, showing all the bifurcations that occur when we vary $D$. 
Similarly, we can obtain the one-parameter bifurcation diagram with $x_1$ or $x_2$ on the $y$-axis.
This diagram corresponds to a vertical line of equation $S_{in}=1$ in the operating diagram in \Cref{FigDOMaple,Fig-DOMatc}.
A magnification of the one-parameter bifurcation diagram is illustrated in \cref{Fig-DB}(b) showing the transcritical bifurcations occur at $\sigma_3$ and $\sigma_4$ which describe the appearance and disappearance of the unique positive steady state $\mathcal{E}^*$.

First, by decreasing $D$, the one-parameter bifurcation diagram in \cref{Fig-DB} illustrates the transcritical bifurcation occurring at $D=\sigma_1=4$ between $\mathcal{E}_0$ and $\mathcal{E}_1$.
Indeed, the washout $\mathcal{E}_0$ loses its stability while the steady state of exclusion of the second species $\mathcal{E}_1$ emerges LES.
Decreasing $D$ further, the steady state of exclusion of the first species $\mathcal{E}_2$ appears unstable at $D=\sigma_2\simeq 1.0666$ via a second transcritical bifurcation with the steady state $\mathcal{E}_0$.
Next, a third transcritical bifurcation occurs at $D=\sigma_3\simeq 0.6447$ that gives birth to an unstable positive steady state $\mathcal{E}^*$ where $\mathcal{E}_2$ changes its stability and becomes LES. 
Finally, decreasing $D$ further, once again a fourth transcritical bifurcation occurs at $D=\sigma_4\simeq 0.352$ between $\mathcal{E}^*$ and $\mathcal{E}_1$ where the steady state $\mathcal{E}^*$ disappears while $\mathcal{E}_1$ remains unstable.
In the following proposition, we summarize the study of the one-parameter bifurcation diagram in $D$ corresponding to the operating diagram in \cref{FigDOMaple,Fig-DOMatc} when the input concentration $S_{in}$ is fixed.
\begin{proposition}
Let $S_{in}=1$. 
Using the biological parameter values provided in \cref{Tab-Allpar} and the specific growth rates \cref{SpeciFunc}, the existence and stability of steady states of \cref{DDInterSpecMod} according to the control parameter $D$ are given in \cref{Tab-ExiStabBD} where the bifurcation values $\sigma_i$, $i=1,\ldots,4$ and the corresponding nature of the bifurcations are defined in \cref{Tab-Sigmai}.
\begin{table}[ht]
\caption{Critical values $\sigma_i$, $i=1,\ldots,4$ of $D$ and the corresponding nature of bifurcations when $S_{in}=1$. \label{Tab-Sigmai}} 
\begin{center}
\begin{tabular}{ @{\hspace{1mm}}l@{\hspace{2mm}} @{\hspace{2mm}}l@{\hspace{2mm}}  @{\hspace{2mm}}l@{\hspace{2mm}}  @{\hspace{1mm}}l@{\hspace{1mm}} }
Definition                 &   Value  &   Bifurcation       \\  
$\sigma_1=f_1(S_{in},0)$   & 4        &    TB    \\ \hline 
$\sigma_2=f_2(S_{in},0)$   & 1.0666   &    TB    \\ \hline 
$\sigma_3$ is the solution of equation $\tilde{x}_2(S_{in},D)=\bar{x}_2(S_{in},D)$  & 0.6447  &  TB     \\ \hline 
$\sigma_4$ is the solution of equation $\tilde{x}_1(S_{in},D)=\bar{x}_1(S_{in},D)$  & 0.352   &  TB     \\ 
\end{tabular}
\end{center}
\end{table}
\begin{table}[ht]
\caption{Existence and stability of all steady states according to $D$ where $\sigma_i$, $i=1,\ldots,4$ are defined in \cref{Tab-Sigmai}. \label{Tab-ExiStabBD}} 
\begin{center}
\begin{tabular}{ @{\hspace{1mm}}l@{\hspace{2mm}} @{\hspace{2mm}}l@{\hspace{2mm}}  @{\hspace{2mm}}l@{\hspace{2mm}}   @{\hspace{2mm}}l@{\hspace{2mm}}  @{\hspace{1mm}}l@{\hspace{1mm}} }
Interval of $D$   & $\mathcal{E}_0$  & $\mathcal{E}_1$ & $\mathcal{E}_2$ &$\mathcal{E}^*$ \\  
$(0,\sigma_4)$         &      U      &      U          &      S          &              \\ \hline
$(\sigma_4,\sigma_3)$  &      U      &      S          &      S          &      U        \\ \hline 
$(\sigma_3,\sigma_2)$  &      U      &        S        &      U          &                \\ \hline 
$(\sigma_2,\sigma_1)$  &      U      &        S        &                 &                 \\ \hline 
$(\sigma_1,+\infty)$   &      S      &                 &                 &                  \\
\end{tabular}
\end{center}
\end{table}
\end{proposition}

\begin{figure}[!ht]
\setlength{\unitlength}{1.0cm}
\begin{center}
\begin{picture}(6.3,6)(0,0)
\put(-4,0){\rotatebox{0}{\includegraphics[width=8cm,height=6cm]{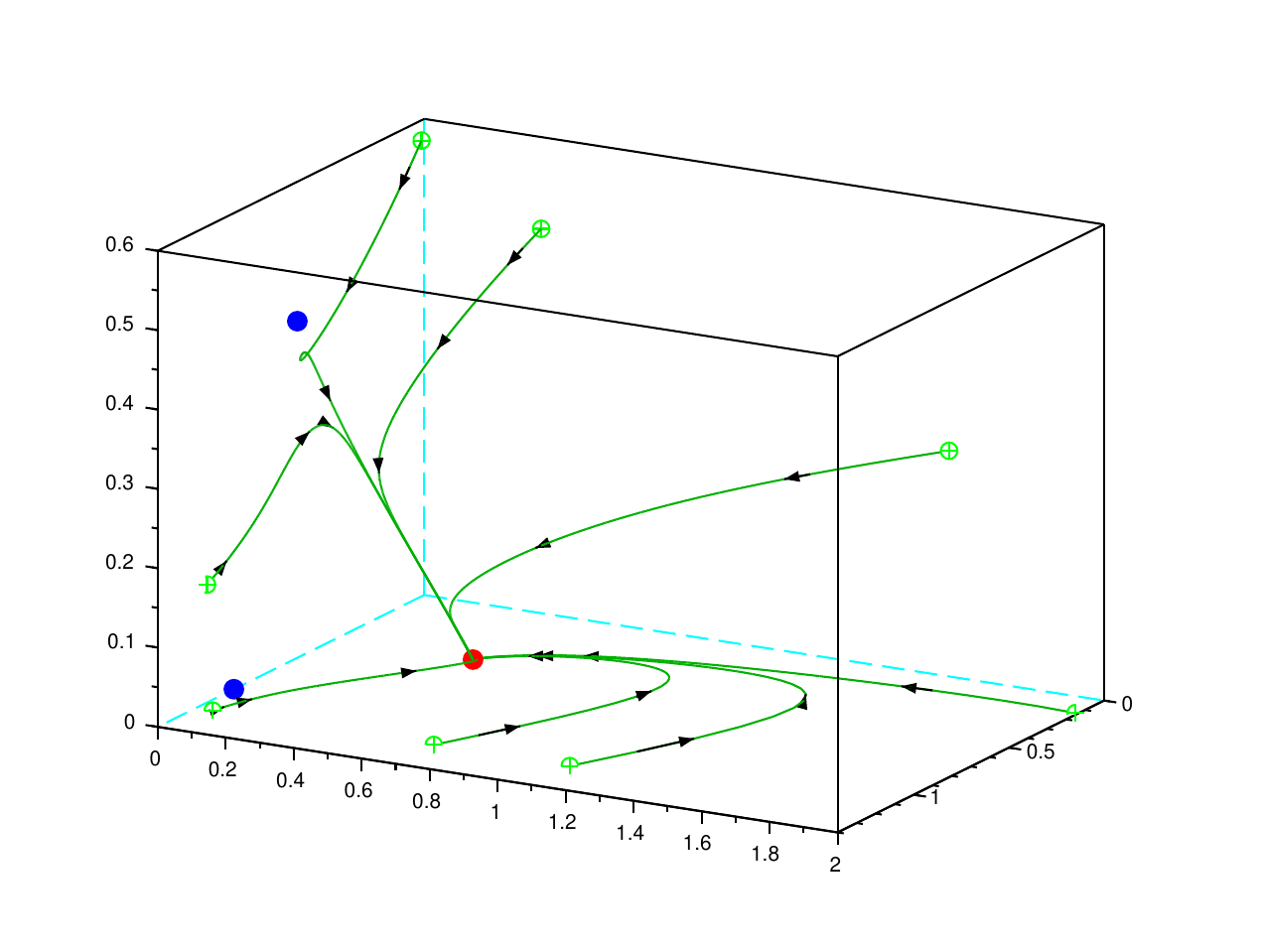}}}
\put(3.1,0){\rotatebox{0}{\includegraphics[width=8cm,height=6cm]{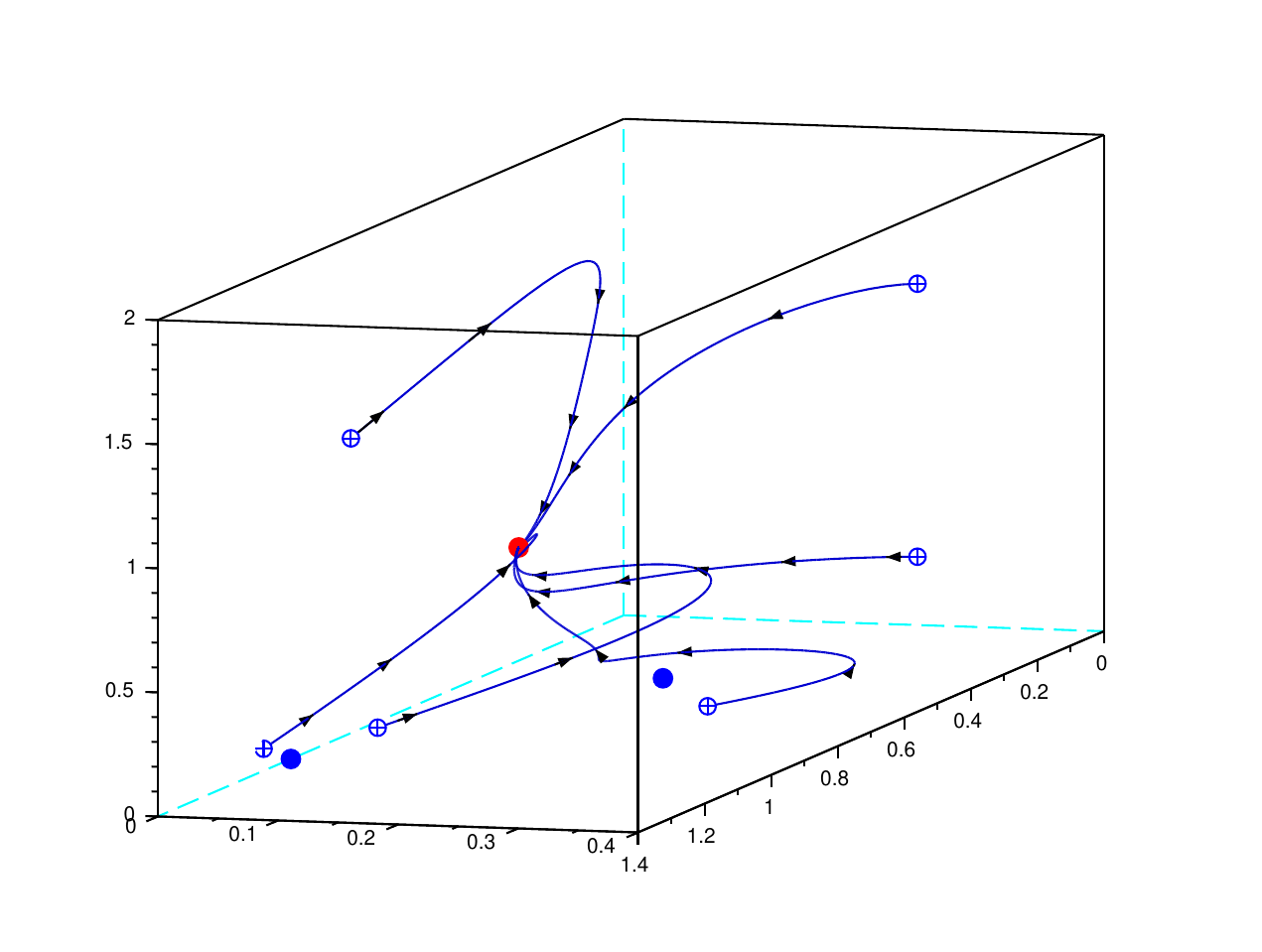}}}
\put(-1.2,-5.65){\rotatebox{0}{\includegraphics[width=9cm,height=6cm]{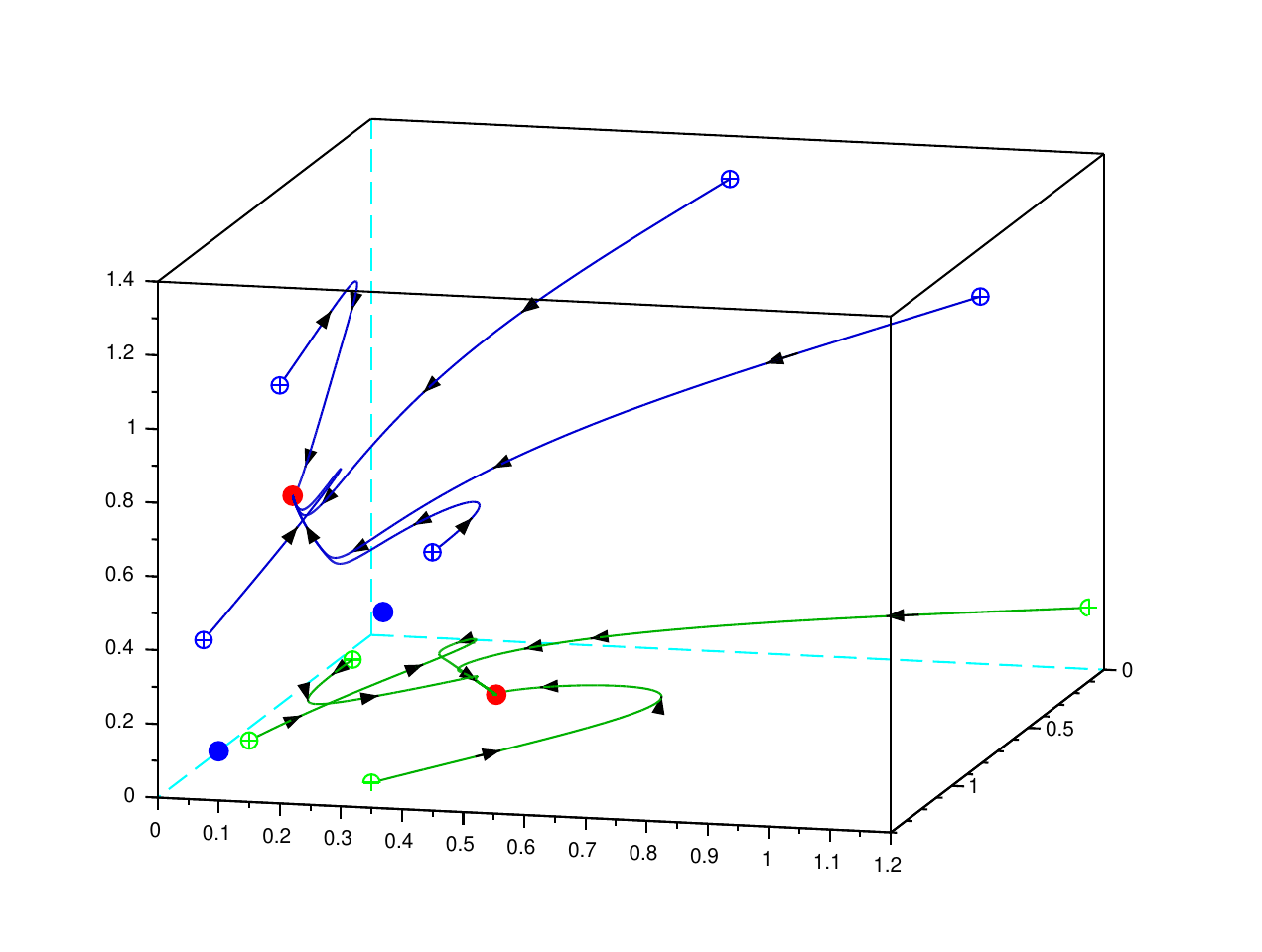}}}
\put(-0.2,5.45){\sc (a)}	
\put(-1.5,0.65){{\sc  $x_1$   }}
\put(-3.65,2.75){{\sc  $x_2$   }}
\put(2.2,0.85){{\sc  $S$   }}
\put(-0.9,2.){{\sc  {\color{red}$\mathcal{E}_1$  } }}
\put(-2.65,1.85){{\sc {\color{blue}$\mathcal{E}_0$  }}}
\put(-2.65,3.9){{\sc {\color{blue}$\mathcal{E}_2$  }}}
\put(6.3,5.45){{\sc (b)}}
\put(5,0.45){{\sc  $x_1$   }}
\put(3.5,2.45){{\sc  $x_2$   }}
\put(9,1.){{\sc  $S$   }}
\put(7.2,1.25){{\sc  {\color{blue}$\mathcal{E}_1$  } }}
\put(5,0.95){{\sc {\color{blue}$\mathcal{E}_0$  }}}
\put(6,2.7){{\sc {\color{red}$\mathcal{E}_2$  }}}
\put(3,-0.2){{\sc (c)}}
\put(2,-5.3){{\sc  $x_1$   }}
\put(-0.8,-3){{\sc  $x_2$   }}
\put(6,-4.7){{\sc  $S$   }}
\put(1.85,-4.25){{\sc  {\color{red}$\mathcal{E}_1$  } }}
\put(0.2,-4.15){{\sc {\color{blue}$\mathcal{E}_0$  }}}
\put(0.4,-2.85){{\sc {\color{red}$\mathcal{E}_2$  }}}
\put(1.2,-3.5){{\sc {\color{blue}$\mathcal{E}^*$  }}}
\end{picture}
\end{center}
\vspace{5cm}
\caption{SCILAB: the three-dimensional space $(S, x_1, x_2)$ of \cref{DDInterSpecMod} according to regions in the operating diagram of \cref{FigDOMaple,Fig-DOMatc} when $S_{in}=1$. (a) Convergence to $\mathcal{E}_1$ for $(S_{in},D)\in \mathcal{J}_5$. (b) Convergence to $\mathcal{E}_2$ for $(S_{in},D)\in \mathcal{J}_4$. (c) Bistability of $\mathcal{E}_1$ and $\mathcal{E}_2$ when $(S_{in},D)\in \mathcal{J}_3$. } \label{Fig-3DSim}
\end{figure}
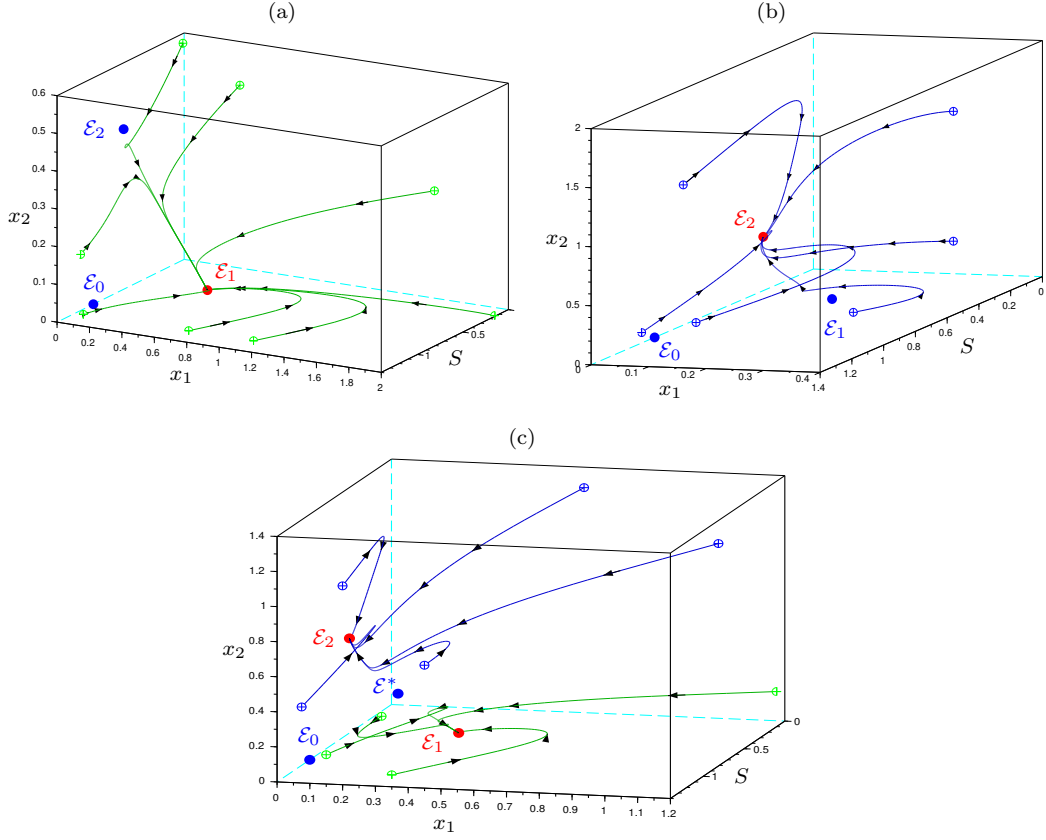

In what follows, we illustrate numerically different cases of stability of steady states using SCILAB \cite{SCILAB2021}.
Using the same set of parameter values in \cref{FigDOMaple,Fig-DOMatc}, the numerical simulations in \cref{Fig-3DSim} illustrate the trajectories in the three-dimensional phase space $(S,x_1,x_2)$.
\Cref{Fig-3DSim}(a) shows that the solutions of system \cref{DDInterSpecMod} converge to $\mathcal{E}_1$ for several positive initial conditions when $(S_{in},D)=(1,0.7)\in \mathcal{J}_5$.
In this case, there exist three steady states given by 
$$
\mathcal{E}_0=(1,0,0),\quad \mathcal{E}_1\simeq( 0.4607,0.4015 ,0),\quad   \mathcal{E}_2\simeq( 0.6666,0,0.4242),
$$
where $\mathcal{E}_1$ is LES while $\mathcal{E}_0$ and $\mathcal{E}_2$ are unstable. The numerical simulations were done for several positive initial conditions permitting us to conjecture the convergence to $\mathcal{E}_1$ for any positive initial condition.
Next, we choose $(S_{in},D)=(1,0.2)\in \mathcal{J}_4$. 
System \cref{DDInterSpecMod} admits the following three steady states 
$$
\mathcal{E}_0=(1,0,0),\quad \mathcal{E}_1\simeq(0.3987,0.1431 ,0),\quad  \mathcal{E}_2\simeq(0.3157,0,0.4561),
$$
where $\mathcal{E}_2$ is LES while $\mathcal{E}_0$ and $\mathcal{E}_1$ are unstable. \Cref{Fig-3DSim}(b) illustrates the convergence towards $\mathcal{E}_2$ for several positive initial conditions.
Finally, when $(S_{in},D)=(1,0.5)$ is chosen in the bistability region $\mathcal{J}_3$, the washout steady state is given by $\mathcal{E}_0=(1,0,0)$ and the other steady states are given by  
$$
\mathcal{E}_1\simeq(0.4354,0.3136 ,0),\quad \mathcal{E}_2\simeq(0.5142,0,0.5396), \quad   \mathcal{E}^*=(0.5181,0.1491,0.2371),
$$
where $\mathcal{E}_0$ and $\mathcal{E}^*$ are unstable while $\mathcal{E}_1$ and $\mathcal{E}_2$ are stable.
In fact, the system exhibits bistability with two basins of attraction, which are separated by the stable manifold of saddle point $\mathcal{E}^*$. 
One basin attracts the solutions to the steady state $\mathcal{E}_1$ and the other to the steady state $\mathcal{E}_2$, as shown in \cref{Fig-3DSim}(c).
The three cases of \cref{Cases123} and their correspondences in the three cases considered in \cref{Fig-3DSim} (with the stable steady states and the regions of the operating diagram) are summarized in \cref{Tab-NumSim}.
\begin{table}[ht]
\caption{The three cases of \cref{Cases123} and their correspondences illustrated in \cref{Fig-3DSim}. \label{Tab-NumSim}} 
\begin{center}
\begin{tabular}{ @{\hspace{1mm}}l@{\hspace{2mm}} @{\hspace{2mm}}l@{\hspace{2mm}}  @{\hspace{2mm}}l@{\hspace{2mm}}  @{\hspace{1mm}}l@{\hspace{1mm}} }
Case              & Region           & Stable steady state  &  Figure     \\ 
$\mbox{Case\:1}$  & $\mathcal{J}_5$  & $\mathcal{E}_1$       &  \cref{Fig-3DSim}(a)                  \\ \hline
$\mbox{Case\:2}$  & $\mathcal{J}_3$  & $\mathcal{E}_1$ and $\mathcal{E}_2$ &  \cref{Fig-3DSim}(c)    \\ \hline  
$\mbox{Case\:3}$  & $\mathcal{J}_4$  & $\mathcal{E}_2$       &  \cref{Fig-3DSim}(b)                
\end{tabular}
\end{center}
\end{table}
\section{Conclusion}	                              \label{Sec-Conc}
In this work, we have extended the intra-specific density-dependent model of \cite{ElhajjiIJB2018} where two species compete for a single resource
by allowing distinct removal rates.
Our study considers the effect of mortality terms, which were neglected in the theoretical analysis of \cite{ElhajjiIJB2018}.
The three-dimensional system cannot be reduced to a two-dimensional one. 

On the other hand, this model was studied in Fekih-Salem et al. \cite{FekihMB2017} in the case of intra- and interspecific density-dependent growth rates. However, the study of the particular case with only interspecific interference between the two populations of microorganisms was not investigated in Fekih-Salem et al. \cite{FekihMB2017}. 
Our study has allowed us to understand whether intraspecific interference is necessary to guarantee persistence between microbial species.


Under general interspecific density-dependent growth functions, we gave an exhaustive study of the behavior of system \cref{DDInterSpecMod}.
More precisely, the necessary and sufficient conditions of existence and stability of all steady states are determined according to the control parameters $S_{in}$ and $D$.
We have shown that system \cref{DDInterSpecMod} admits at most four types of steady states: the washout steady state $\mathcal{E}_0$ which always exists, two boundary steady states $\mathcal{E}_1$ and $\mathcal{E}_2$ in which one of species degenerates and the other persists and a coexistence steady state of two species $\mathcal{E}^*$.
We have proved that whenever $\mathcal{E}^*$ exists, it is unique and unstable.
In fact, the existence and the uniqueness of the positive steady state $\mathcal{E}^*$ depends on the relative position of the values $\tilde{x}_i$ and $\bar{x}_i$, $i=1,2$ (as shown in \cref{Fig-CesesPos}).
In addition, a positive steady state exists if and only if both boundary steady states are LES.
In this case, the model exhibits the bi-stability. 

Our result on the instability of the positive steady state in this paper can be deduced from those in \cite{FekihMB2017}. 
More precisely, \cite[condition (28)]{FekihMB2017} holds in the particular case without intraspecific interference and \cite[Theorem 1]{FekihMB2017} implies in this case that a positive steady state, if it exists, is always unique and unstable.

Using our mathematical analysis of model \cref{DDInterSpecMod}, we have determined, first, theoretically the operating diagram to provide a complete description of the asymptotic behavior of the process with respect to the control parameters $S_{in}$ and $D$.
Our description of the operating diagram in the $(S_{in},D)$-plane using the specific growth functions \cref{SpeciFunc} shows that this diagram is divided into six regions: $\mathcal{J}_0$ corresponds to the washout of two species; the regions $\mathcal{J}_1$ and $\mathcal{J}_5$ [resp. $\mathcal{J}_2$ and $\mathcal{J}_4$] correspond to the competitive exclusion of the second [resp. first] species; the bi-stability region $\mathcal{J}_3$ of the boundary steady states with the existence of a unique unstable positive steady state.

On the other hand, we have established numerically the operating diagram of \cref{DDInterSpecMod} using the software MatCont \cite{MATCONT2023}.
We found the same operating diagram obtained theoretically from the existence and stability conditions.
However, this numerical method was able to detect new types of two-parameter bifurcations: BT, which corresponds to a Bogdanov-Takens bifurcation, and ZH, which corresponds to a Zero-Hopf bifurcation (see \cref{Fig-DOMatc}).

When only the dilution rate $D$ varies, the study of the one-parameter bifurcation diagram has illustrated the nature of bifurcations of all steady states by passing through the various regions of the operating diagram in \cref{FigDOMaple,Fig-DOMatc}.
It was shown that all steady states can appear or disappear only through transcritical bifurcations.
The numerical simulations have shown in three-dimensional phase space $(S,x_1,x_2)$ that the system can exhibit bi-stability with convergence either to the exclusion of the first species or to the exclusion of the second species.

Finally, our mathematical analysis proves that the outcome of competition satisfies the CEP which predicts that only one species can exist. 
Thus, our study allows us to say that the addition of mortality terms of the species in the interspecific density-dependent model with a mutual inhibitory relationship, could not lead to the coexistence around a stable positive steady state or stable limit cycle.
Comparing our results with those in \cite{FekihMB2017}, we can conclude that interspecific interference is not sufficient to guarantee coexistence between species, contrary to intraspecific interference.
\appendix
\section{Proofs}  \label{AppendixA}
\begin{proof}[\bf{Proof of \cref{Lemma-EquivXiFi}}]
Under \cref{hyp2}, it is easy to see that the function $x_1\longmapsto f_2(S_{in}-D_1x_1/D,x_1)-D_2$ is decreasing. From the definition of $\bar{x}_1$ in \cref{Tab-BECAN}, we have
$$
\bar{x}_1<\tilde{x}_1 \quad \Leftrightarrow \quad f_2(S_{in}-D_1\tilde{x}_1/D,\tilde{x}_1)-D_2 < f_2(S_{in}-D_1\bar{x}_1/D,\bar{x}_1)-D_2 = 0.
$$
From the definition of $\tilde{x}_1$ in \cref{Tab-BECAN}, we obtain the first equivalence in \cref{EquivXiFi}.
Similar arguments apply to the second equivalence in \cref{EquivXiFi}.
\end{proof}
\begin{proof}[\bf{Proof of \cref{Prop-Exis}}]
\begin{itemize}[leftmargin=*]
\item When $x_1=x_2=0$, from the first equation of \cref{SystSS}, one has $S=S_{in}$.
So, $\mathcal{E}_0=(S_{in},0,0)$, it always exists.
\item For the boundary steady state $\mathcal{E}_i$, $x_i>0$ and $x_j=0$, $i=1,2$, $j=1,2$, $j\neq i$.
From system \cref{SystSS}, the components $\tilde{S}_i$ and $\tilde{x}_i$ are the solutions of 
\begin{align}
D\left(S_{in}-\tilde{S}_i\right) &= D_i\tilde{x}_i  \label{EquSixitilde}\\
f_i(\tilde{S}_i,0)    &= D_i. \label{EqulambdiBEC}
\end{align}
Since the function $S\longmapsto f_i(S,0)$ is increasing (under \cref{hyp2}), it follows from equation \cref{EqulambdiBEC} and \cref{Tab-BECAN} that
$$
\tilde{S}_i=\lambda_i(D).
$$
Replacing $\tilde{S}_i$ by this expression in \cref{EquSixitilde}, we can see that 
$$
\tilde{x}_i=D(S_{in}-\lambda_i(D))/D_i,
$$ 
which is positive if and only if $S_{in} > \lambda_i(D)$. If $\mathcal{E}_i$ exists, then it is unique.
\item For the positive steady state $\mathcal{E}^*$, $x_1>0$ and $x_2>0$.
From system \cref{SystSS}, the components $S^*$, $x_1^*$ and $x_2^*$ are  the solutions of the set of equations
\begin{align}
D(S_{in}-S)  &= D_1x_1+D_2x_2,    \label{EquSolSet} \\
f_1(S,x_2)   &= D_1,              \label{Equf1=D1}     \\
f_2(S,x_1)   &= D_2.              \label{Equf2=D2}
\end{align}
Hence, from \cref{EquSolSet}, we obtain
\begin{equation}                   \label{ExpSetoil}
S^*= S_{in}-D_1x_1^*/D-D_2x_2^*/D.
\end{equation}
Replacing expression \cref{ExpSetoil} in \cref{Equf1=D1}-\cref{Equf2=D2}, one has $(x_1=x_1^*,x_2=x_2^*)$ is a solution of 
\begin{equation}                     \label{EqExixietoil}
\left\{
\begin{array}{ll}
f_1(S_{in}-D_1x_1/D-D_2x_2/D,x_2)  = D_1, \\
f_2(S_{in}-D_1x_1/D-D_2x_2/D,x_1)  = D_2.
\end{array}
\right.
\end{equation}
From \cref{ExpSetoil}, one see that $S^*$ is positive if and only if $D_1x_1^*/D+D_2x_2^*/D<S_{ in}$.
Thus, \cref{EqExixietoil} has a positive solution in the interior of the set $M$ defined by \cref{SetM}, that is, the curves $\gamma_1$ of $F_1$ and $\gamma_2$ of $F_2$, respectively, has a positive intersection in the interior of $M$ (from \cref{Lem-Fi}).
More precisely, $(x_1=x_1^*,x_2=x_2^*)$ is a solution of equation
$$
x_2=F_1(x_1)\quad \mbox{and} \quad x_2=F_2(x_1).
$$
Let $F(x_1)$ be the function defined by
\begin{equation}                                          \label{DefF}
F(x_1):=F_1(x_1)-F_2(x_1).
\end{equation}
Using the derivatives \cref{InegF1prim} of the function $F_1$ and \cref{InegF2prim} of $F_2$, respectively, a straightforward calculation shows that
$$
F'_1(x_1)-F'_2(x_1)=\frac{D(DGH+D_1FG+D_2EH)}{D_2F(D_2E+DG)}>0, \quad \mbox{for all} \quad x_1<\min (\tilde{x}_1,\bar{x}_1).
$$
Therefore, the function $F$ is increasing from $F(0)=\bar{x}_2-\tilde{x}_2$ to $F(\bar{x}_1)=F_1(\bar{x}_1)-F_1(\tilde{x}_1)$ since $F_2(\bar{x}_1)=0=F_1(\tilde{x}_1)$. 
Hence, $F(\bar{x}_1)>0$ if and only if $\bar{x}_1<\tilde{x}_1$.
Consequently, there exists an $x_1=x_1^*\in(0,\bar{x}_1)$ satisfying \cref{DefF} if and only if
$$
\bar{x}_1<\tilde{x}_1 \quad \mbox{and} \quad  \bar{x}_2<\tilde{x}_2.
$$
If $\tilde{x}_1$ exists, then it is unique. This completes the proof.
\end{itemize}
\end{proof}
\begin{proof}[\bf{Proof of \cref{Prop-Stb}}]
Here, we study the stability conditions of all steady states of \cref{DDInterSpecMod}. The Jacobian matrix of \cref{DDInterSpecMod} at $(S,x_1,x_2)$ takes the form
\begin{equation*}
J=
\begin{pmatrix}
-D-x_1E-x_2F    &      -f_1(S,x_2)+x_2H   &   x_1G-f_2(S,x_1) \\
x_1E            &     f_1(S,x_2)-D_1      &   -x_1G \\
x_2F            &     -x_2H               &    f_2(S,x_1)- D_2 \\
\end{pmatrix}
\end{equation*}
where $E$, $F$, $G$ and $H$ are given by \cref{ExprEFGH}.
\begin{itemize}[leftmargin=*]
\item At $\mathcal{E}_0=(S_{ in},0,0)$, the characteristic polynomial is
$$
P(\lambda)=(-D-\lambda)(f_1(S_{in},0)-D_1-\lambda)(f_2(S_{in},0)-D_2-\lambda).
$$
Therefore, $\mathcal{E}_0$ is LES if and only if $f_i(S_{in},0)<D_i$, for $i=1,2$.
Using the definition of $\lambda_i(D)$ in \cref{Tab-BECAN}, we conclude that $\mathcal{E}_0$ is LES if and only if $S_{in}<\lambda_i(D)$, for $i=1,2$.
\item At $\mathcal{E}_1=(\lambda_1(D),D(S_{in}-\lambda_1(D))/D_1,0)$, the characteristic polynomial is
$$
P(\lambda)=(f_2(\lambda_1(D),D(S_{in}-\lambda_1(D))/D_1)-D_2-\lambda)(\lambda^2+\beta_1\lambda+\beta_2),
$$
where $\beta_1=D(1+E(S_{in}-\lambda_1(D))/D_1)$ and $\beta_2=DE(S_{in}-\lambda_1(D))$.
Since $\beta_1>0$ and $\beta_2>0$, the real parts of the roots of the quadratic factor are negative.
Therefore, $\mathcal{E}_1$ is LES if and only if $f_2(\lambda_1(D),D(S_{in}-\lambda_1(D))/D_1)<D_2$. 
\item Similarly, at $\mathcal{E}_2=(\lambda_2(D),0,D(S_{in}-\lambda_2(D))/D_2)$, the characteristic polynomial is
$$
P(\lambda)=(f_1(\lambda_2(D),D(S_{in}-\lambda_2(D))/D_2)-D_1-\lambda)(\lambda^2+\delta_1\lambda+\delta_2),
$$
where $\delta_1=D(1+F(S_{in}-\lambda_2(D))/D_2)$ and $\delta_2=DF(S_{in}-\lambda_2(D))$.
Since $\delta_1>0$ and $\delta_2>0$, the real parts of the roots of the quadratic factor are negative.
Therefore, $\mathcal{E}_2$ is LES if and only if $f_1(\lambda_2(D),D(S_{in}-\lambda_2(D))/D_2)<D_1$.
\item The characteristic polynomial of the Jacobian matrix $J$ evaluated at $\mathcal{E}^*$ is given by:
$$
P(\lambda)=\lambda^3+c_1\lambda^2+c_2\lambda+c_3,
$$
where 
\begin{align*}                          
\begin{split}
& c_1 = D+Ex_1^*+Fx_2^*,\quad c_2=D_1Ex_1^*+D_2Fx_2^*-(GH+FG+EH)x_1^*x_2^*,\\
& c_3 = -(DGH+D_1FG+D_2EH)x_1^*x_2^*.
\end{split}
\end{align*}
Hence, the coefficient $c_3<0$. According to the Routh--Hurwitz criterion, $\mathcal{E}^*$ is always unstable, whenever it exists.
\end{itemize}
\end{proof}
\section{Parameters used in numerical simulations} \label{AppendixB}
All the values of the parameters used in the numerical simulations are provided in \cref{Tab-Allpar}.
For the numerical simulations, we have used MAPLE \cite{MAPLE18} to plot \cref{Fig-CesesPos,FigDOMaple,Fig-DB}, MATCONT \cite{MATCONT2023} for \cref{Fig-DOMatc} and SCILAB \cite{SCILAB2021} for \cref{Fig-3DSim}.
\begin{table}[!ht]
\caption{Parameter values used for model \cref{DDInterSpecMod} when the growth rates $f_1$ and $f_2$ are given by \cref{SpeciFunc}.}   \label{Tab-Allpar}
\vspace{-0.2cm}
\centering{\sl
\begin{tabular}{ @{\hspace{1mm}}l@{\hspace{2mm}} @{\hspace{2mm}}l@{\hspace{2mm}} @{\hspace{2mm}}l@{\hspace{2mm}} @{\hspace{2mm}}l@{\hspace{2mm}}
                 @{\hspace{2mm}}l@{\hspace{2mm}} @{\hspace{2mm}}l@{\hspace{2mm}} @{\hspace{2mm}}l@{\hspace{2mm}} @{\hspace{2mm}}l@{\hspace{2mm}}
                 @{\hspace{2mm}}l@{\hspace{2mm}} @{\hspace{2mm}}l@{\hspace{2mm}} @{\hspace{2mm}}l@{\hspace{2mm}} @{\hspace{2mm}}l@{\hspace{2mm}} @{\hspace{2mm}}l@{\hspace{2mm}}}
Parameter     & $ m_1$ & $K_1$ & $\beta_1$ & $m_2$  & $K_2$  &  $\beta_2$  &  $\alpha_1$  &  $\alpha_2$  &  $a_1$   &  $a_2$  & $S_{in}$  & $D$ \\ \hline

\begin{tabular}{@{\hspace{0mm}}l@{\hspace{0mm}}}
Figs. \ref{Fig-CesesPos}(a), \ref{Fig-3DSim}(a)   \\
Figs. \ref{Fig-CesesPos}(b), \ref{Fig-3DSim}(c) \\
Figs. \ref{Fig-CesesPos}(c), \ref{Fig-3DSim}(b) \\
Figs. \ref{FigDOMaple}, \ref{Fig-DOMatc} \\
Fig.  \ref{Fig-DB}
\end{tabular}
              & 4 & 1.5 & 1.2 & 2.2 & 2 & 0.1 & 0.2 & 0.5 & 0.8 & 0.2 &
\begin{tabular}{@{\hspace{0mm}}l@{\hspace{0mm}}}
1     \\
1    \\
1    \\
var  \\
1
\end{tabular}
&
\begin{tabular}{@{\hspace{0mm}}l@{\hspace{0mm}}}
0.7     \\
0.5    \\
0.2    \\
  var \\
var
\end{tabular}
\end{tabular}}
\end{table}



\end{document}

%% file: MF_shared.tex
\usepackage{lipsum}
\usepackage{amsfonts}
\usepackage{graphicx}
\usepackage{epstopdf}
\usepackage{algorithmic}
\usepackage{xcolor}

\def\sc{\scriptsize}
\def\sl{\small}
\def\ds{\displaystyle}

\ifpdf
  \DeclareGraphicsExtensions{.eps,.pdf,.png,.jpg}
\else
  \DeclareGraphicsExtensions{.eps}
\fi

\graphicspath{{Images/}}

\newcommand{{\sSN}}{\ensuremath{ {\mbox{\tiny{SN}}} }}

\usepackage{enumitem}
\setlist[enumerate]{leftmargin=.5in}
\setlist[itemize]{leftmargin=.5in}

\newsiamremark{remark}{Remark}
\newsiamremark{hypothesis}{Hypothesis}
\crefname{hypothesis}{Hypothesis}{Hypotheses}
\newsiamthm{claim}{Claim}

\headers{Interspecific density-dependent model}{T. MTAR, R. FEKIH-SALEM}

\title{Analysis and operating diagram of an interspecific density-dependent model 
\thanks{Submitted to the editors 2024-01-05.
\funding{This work was supported by the Euro-Mediterranean research network TREASURE (\url{http://www.inra.fr/treasure}).}
}}

\author{
     TAHANI MTAR  \thanks{University of Tunis El Manar, National Engineering School of Tunis, LAMSIN, Tunisia (\email{tahani.mtar@enit.utm.tn}).}
\and RADHOUANE FEKIH-SALEM \footnotemark[1]~\footnotemark[2]~\thanks{University of Monastir, Higher Institute of Computer Science of Mahdia, Tunisia     (\email{radhouane.fekih-salem@enit.utm.tn}).}
}
\usepackage{amsopn}
